\documentclass{article}
\usepackage{arxiv}
\usepackage[utf8]{inputenc} % allow utf-8 input
\usepackage[T1]{fontenc}    % use 8-bit T1 fonts
\usepackage{hyperref}       % hyperlinks
\usepackage{url}            % simple URL typesetting
\usepackage{booktabs}       % professional-quality tables
\usepackage{amsfonts}       % blackboard math symbols
\usepackage{nicefrac}       % compact symbols for 1/2, etc.
\usepackage{microtype}      % microtypography
\usepackage{lipsum}		% Can be removed after putting your text content
\usepackage{graphicx}
\usepackage{natbib}
\usepackage{doi}
\usepackage{tabularx}
\usepackage{multicol}
\usepackage{multirow}
\usepackage{pdflscape}
\usepackage{array}
\usepackage{tikz}
\usepackage{tabularx}
\usepackage{multicol}
\usepackage{multirow}
\usepackage{pdflscape}
\usepackage{array}
\usepackage{threeparttable}
\usepackage{tikz}
\usetikzlibrary{arrows, shapes, positioning, calc, shadows}
\usepackage{ragged2e}
\usepackage{amsmath}
\usepackage{subcaption}
\usepackage{caption}

\usepackage{microtype}

\usepackage[para]{footmisc}
\usepackage{setspace}
\usepackage[format=hang,singlelinecheck=false]{caption}
\usepackage{ragged2e}

\usepackage{longtable}

\title{A Two-Stage Stochastic  Model for Road-Rail Intermodal Freight Transportation Under Demand and Capacity Uncertainty}

\author{ \href{https://orcid.org/0009-0007-7093-748X}{\includegraphics[scale=0.06]{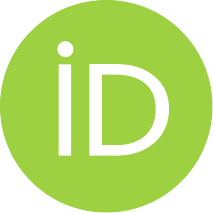}\hspace{1mm}Jeremiah ~Gbadegoye}\\
	Department of Industrial and Systems Engineering\\
	University of Tennessee\\
	Knoxville, TN 37996 \\
	\texttt{jgbadego@vols.utk.edu} \\
    \And
	\href{https://orcid.org/0000-0001-7465-7783}{\includegraphics[scale=0.06]{orcid.pdf}\hspace{1mm}Mustafa C.~Camur} \\
	Department of Industrial and Systems Engineering\\
	University of Tennessee\\
	Knoxville, TN 37996 \\
	\texttt{mcamur@utk.edu} \\
	\And
	\href{https://orcid.org/0000-0003-1990-0159}{\includegraphics[scale=0.06]{orcid.pdf}\hspace{1mm}Xueping ~Li}\\
	Department of Industrial and Systems Engineering\\
	University of Tennessee\\
	Knoxville, TN 37996 \\
	\texttt{Xueping.Li@utk.edu} \\
}

\renewcommand{\shorttitle}{Two-Stage Stochastic  Model for Road-Rail Intermodal Freight Transportation}

\hypersetup{
pdftitle={A template for the arxiv style},
pdfsubject={q-bio.NC, q-bio.QM},
pdfauthor={David S.~Hippocampus, Elias D.~Striatum},
pdfkeywords={First keyword, Second keyword, More},
}

\begin{document}
\maketitle

\begin{abstract}
With the steady increase in global logistics and freight transport demand, the need for efficient and sustainable intermodal transport systems becomes increasingly important. This study addresses the optimization of container movement by intermodal transport with fixed train schedules. We emphasize the integration of road-rail intermodal transport amid uncertain demand and train (spot) capacities. A two-stage stochastic optimization model is developed to strategically manage the transportation of containers from multiple origins to designated intermodal hubs. By leveraging spot capacities at train stations and addressing uncertainties in demand and train capacity, the model integrates Conditional Value-at-Risk (CVaR) to balance cost efficiency and risk management, enabling robust decision-making under uncertainty. The model’s objectives encompass minimizing transportation costs, mitigating carbon emissions, and enhancing the reliability of containerized freight movement across the network. A comprehensive case study using real-world data demonstrates the practical applicability of the model, highlighting its effectiveness in reducing operational costs, minimizing environmental impacts, and providing actionable insights for stakeholders to navigate the trade-offs between expected costs and risk management in dynamic intermodal transport settings.
 
\end{abstract}

% keywords can be removed
\keywords{Intermodal Freight Transportation \and Stochastic Optimization \and Decarbonization  \and Road-Rail Transport}

\section{Introduction}

Intermodal freight transportation involves the movement of freight using at least two transportation modes without changing the load unit throughout the process. This approach leverages the cost-effectiveness and high capacity of modes like rail or sea for long-haul segments while relying on trucks for shorter distances. It typically requires consolidating multiple shipments into containers or single units at facilities such as warehouses or manufacturing points. These consolidated units are then transported to intermodal hubs for shipping and eventually transferred to final intermodal nodes for deconsolidation and final delivery. The scale and importance of intermodal freight are evident in the fact that, in 2021, the United States (US) alone transported an average of 53.6 million tons of freight daily, valued at over \$54 billion \citep{bts2021movinggoods}, emphasizing it's significant role in the national economy.

Despite its advantages, intermodal freight transport faces several challenges due to the involvement of multiple stakeholders and decision-makers, each with their own priorities and objectives. In the context of our work, a  significant challenge is the modal shift from road to rail transport \citep{el2022logistic}. Integrating these modes efficiently requires overcoming additional handling and time requirements compared to direct road transport \citep{gorccun2024integrated}. For longer distances, rail's high fixed costs are offset by increased efficiency, while road transport becomes less cost-effective as distance increases \citep{zgonc2019impact}.

Synchronized operations within intermodal networks are very crucial, especially when integrating modes with fixed schedules, such as rail transport with more flexible modes like road transport. For instance, trucks must arrive at rail terminals in time to accommodate service time for modal shift and align with train departure schedules. Achieving this requires precise planning and real-time adjustments to account for potential disruptions such as traffic congestion, accidents, or weather-related delays. Effective coordination minimizes delays and ensures the overall efficiency of the intermodal transport system.

Another critical challenge is the uncertainty in market demand and spot capacities at intermodal hubs. Demand fluctuations can arise from various factors, including booking changes, cancellations, unforeseen events, and seasonal variations \citep{wang2021optimizing}. Compounding this uncertainty is the stochastic nature of spot capacities, which depend not only on competing cargo from other companies but also on the variability of loading and unloading activities at each station \citep{sun2016holding}. These dynamics create a complex environment where maintaining operational efficiency requires careful planning and adaptability.

A key aspect influencing spot capacity allocation is the role of long-term contracts, which are agreements between transport providers and large freight forwarders, ensuring a guaranteed capacity over a time period. Such contracts offer stable and predictable volumes for the transport providers while securing regular service for the customers. However, any unused capacity are often redirected to the spot market. This cascading allocation process introduces uncertainty for spot customers, as the availability of capacity depends on how much of the allocated slots long-term contract holders utilize \citep{wang2021optimizing}. For instance, if contract holders underutilize their slots, the surplus is released for spot demand, but this availability may still be constrained by the need to accommodate delayed or overbooked cargo from other services or other long-term contractors. This interplay between long-term contracts and spot demand underscores the need for flexible and adaptive strategies to ensure optimal utilization of available capacity while addressing the inherent uncertainty in both supply and demand \citep{lu2020fuzzy, delbart2021uncertainty}.

While previous studies have addressed aspects of uncertainty management, facility location, or sustainable logistics separately \citep{guo2020dynamic, hernandez2024evaluation}, an integrated framework that jointly considers stochastic demand, spot capacity fluctuations, and carbon emission reduction remains limited. This gap is further discussed in Section~\ref{Literature review}. Moreover, many existing models assume perfect information or static decisions, overlooking the need for adaptive strategies in real-world intermodal operations.

In this study, we develop a two-stage stochastic optimization framework that explicitly incorporates uncertainty in both demand and train capacities while embedding sustainability considerations into intermodal freight planning. Our model leverages a Conditional Value-at-Risk (CVaR) approach to manage the risks associated with extreme cost outcomes, enabling decision-makers to balance operational costs with risk aversion preferences. In addition, the model penalizes excess carbon emissions, aligning the optimization objectives with broader environmental goals.

To demonstrate the practical relevance of the proposed framework, we apply it to a case study involving major freight rail corridors in the United States. Using real-world data on train schedules, intermodal hub locations, and supply facility distributions, the case study evaluates the model’s performance under multiple uncertainty scenarios and explores the trade-offs between cost, risk, and sustainability objectives.

The remaining sections of this paper are organized as follows. Section \ref{Literature review} presents a literature review on intermodal transport resilience. Section \ref{Problem Description} provides a detailed problem description. Section \ref{Optimization Model} discusses our proposed optimization model. Section \ref{Computational Experiments} presents our computational experiments and insights from sensitivity analysis. Section \ref{Conclusion} concludes the paper and suggests directions for future research.

\section{Literature review} \label{Literature review}

Intermodal freight transport has been widely studied for its potential to improve logistics efficiency, sustainability, and resilience. Recent advancements in optimization techniques, data-driven decision-making, and sustainable transport policies have addressed key challenges such as managing uncertainties, coordinating operations across modes, and minimizing emissions. This section highlights significant recent contributions in these areas, examines emerging approaches, and identifies research gaps, providing context for our proposed model's contributions.

\subsection{Managing Uncertainties in Intermodal Transport}

Managing uncertainties in intermodal transport has been a critical in improving logistics efficiency. \citet{guo2020dynamic} tackled demand fluctuations using a multistage stochastic programming model with a rolling horizon approach, outperforming traditional methods in cost reduction and allocation performance. Likewise, \citet{sun2020fuzzy} applied a fuzzy multi-objective routing model for hazardous materials, using triangular fuzzy numbers and soft time windows to balance risk, cost, and service reliability.

Facility location and disruption management have been used to address transport stochasticity. \citet{badyal2023two} proposed a two-stage stochastic model for terminal location and freight distribution under disruption risks, emphasizing the need of robust planning. \citet{bagheri2024optimizing} compared deterministic, stochastic, and robust models for dry port placement, highlighting the flexibility and environmental benefits of robust optimization under dynamic conditions.

Routing and network design remain pivotal in handling operational variability. \citet{lu2020fuzzy} introduced fuzzy programming to address uncertainties in road-rail transport, improving performance through defuzzification and sensitivity analysis. \citet{yang2024path} developed a bi-level genetic algorithm for multimodal routing, effectively minimizing costs and meeting time-sensitive delivery constraints. In addition, \citet{wang2024research} applied adaptive differential evolution to marine transport planning, significantly cutting delays and cost in multimodal networks.

Finally, decision-making frameworks offer a strategic lens for managing systemic uncertainty. \citet{gandhi2024evaluation} used a hybrid multi-criteria framework to assess modal shift risks and enablers, underscoring the role of policy and partnerships. Similarly, \citet{raad2024hybrid} proposed a fuzzy stochastic model for closed-loop dry port networks, boosting resilience through robust, scenario-based planning. These contributions underscore the growing importance of integrating analytical tools and long-term strategies to build agile, sustainable intermodal systems.

\subsection{Operational Synchronization in Road-Rail Intermodal Transport}

Effective synchronization between road and rail is vital for reducing delays and enhancing intermodal efficiency. \citet{newman2000scheduling} pioneered train schedule optimization using decomposition techniques to streamline operations and meet delivery timelines. Building on this, \citet{liu2024flexible} tackled port-side delays with a flexible yard crane scheduling model, optimizing task sequencing and significantly reducing truck congestion for high rail-volume.

Advances in scheduling and routing have improved network performance. \citet{chen2016schedule} addressed uncertainties in vehicle dispatching at terminals, proving better logistics coordination and fewer delays. \citet{chen2022integrated} developed a bi-level model integrating route planning and timetable adjustment, showing higher capacity use and lower costs. For emergency scenarios, \citet{tong2024optimizing} introduced a hybrid flow shop model with a memetic algorithm to minimize congestion and transshipment delays, highlighting the value of adaptive scheduling in multimodal systems.

Strategic planning and infrastructure upgrades also play a central role. \citet{delgado2021intermodal} proposed a decentralized terminal management framework, lowering costs and improving planning in Portugal’s rail-road network. \citet{li2022multimodal} used a network equilibrium model to identify bottlenecks on China-Europe routes, guiding infrastructure expansion. In contrast, \citet{kumar2020evaluating} highlighted how underdeveloped infrastructure and limited transshipment technologies constrain emerging economies, advocating for targeted investments and policy support.

Finally, addressing systemic inefficiencies is key to sustained coordination. \citet{ke2022managing} proposed robust optimization strategies for hazardous material transport, including rerouting, yard repairs, and third-party carriers, boosting network reliability. \citet{ahmady2022optimizing} added a multi-objective model balancing costs, efficiency, and infrastructure criticality advancing the broader goal of synchronized, resilient intermodal operations.

\subsection{Decarbonization in Intermodal Freight Transport}

Decarbonizing intermodal freight transport has driven a wave of innovation through optimization models, technology, and policy. \citet{cheng2019multimodal} introduced a multimodal path optimization model that minimized emissions and costs by integrating carbon taxation and trading into time-constrained freight systems. Their genetic algorithm-based solution proved effective in maintaining timeliness while reducing carbon output. Likewise, \citet{sun2023modeling} developed a fuzzy nonlinear programming model to optimize carbon-efficient road-rail routing, balancing environmental and operational goals under uncertainty.

Technological innovation and alternative fuels are central to emission reduction. \citet{hernandez2024evaluation} assessed biodiesel, e-fuels, hydrogen, and electric locomotives for rail decarbonization, finding electric systems could cut emissions by 46\% under favorable costs. \citet{matsuyama2024scenario} showed that, in Japan’s hinterland transport, emerging technologies could achieve up to 87\% emission reductions, surpassing the impact of modal shift policies alone.

Policy-driven modal shifts also show significant promise. In Brazil, \citet{nassar2023system} used a System Dynamics model to stress that infrastructure and policy must advance in tandem for meaningful decarbonization. These studies reinforce the role of integrated policy frameworks in supporting sustainable freight. Finally, consolidation and network optimization contribute meaningfully to decarbonization. \citet{goodarzi2024evaluating} proposed a two-stage stochastic model accounting for disruption scenarios, showing strong environmental and cost benefits. By enhancing network reliability and mitigating disruptions, their work highlighted the value of resilience-focused strategies.

\subsection{Our Contribution}

Intermodal freight transport, particularly in road-rail systems, faces significant challenges due to fluctuating demand, constrained spot capacities, and fixed train schedules. To address these complexities, we propose a two-stage stochastic optimization model that supports robust first-stage container preparation decisions, followed by adaptive second-stage recourse actions. This structure offers the flexibility needed to operate efficiently in dynamic, uncertain environments \citep{tolooie2020two, huang2021two}.

Our work advances the literature by unifying key priorities, uncertainty management, operational synchronization, and decarbonization into a single integrated framework. Unlike prior studies that examine these elements in isolation \citep{guo2020dynamic, sun2023modeling}, our model also leverages CVaR to account for extreme cost scenarios, striking a deliberate balance between risk resilience, cost efficiency, and environmental impact.

A central innovation of our approach lies in its dynamic adaptation to real-time disruptions, such as variable train capacities and demand shifts. By incorporating spot market capacities at intermodal hubs, our model enables responsive container allocation and preparation within tight operational windows. This marks a departure from static planning models \citep{badyal2023two}, allowing decisions to adjust fluidly with on-the-ground realities.

Equally important is our explicit integration of emissions constraints. While previous studies explore low-carbon strategies \citep{hernandez2024evaluation}, our model embeds carbon penalties directly into the objective function, ensuring sustainability is not an afterthought but a driving force. This aligns with regulatory imperatives like the Clean Air Act \citep{epa2024cleanair} and California’s aggressive emissions targets \citep{wri2024usstate}, reinforcing the relevance of our work to real-world policy goals.

By combining operational adaptability, environmental measures, and risk-aware planning into a unified stochastic framework, we offer a practical, forward-looking tool for designing resilient and sustainable intermodal systems. The following section outlines the operational context that grounds and motivates our model development.

\section{Problem Description}\label{Problem Description}

We address the optimization challenges faced by transportation and manufacturing firms managing intermodal logistics under fixed train schedules. These firms must transport containers from multiple supply points such as warehouses or plants to intermodal hubs in time for scheduled train departures. A critical element of this process is the use of spot capacities on trains, which are often uncertain due to fluctuating demand and competing allocations (e.g., long-term contracts, operational constraints) \citep{wang2021optimizing}.

\begin{figure}[ht]
\centering
\includegraphics[width=\textwidth]{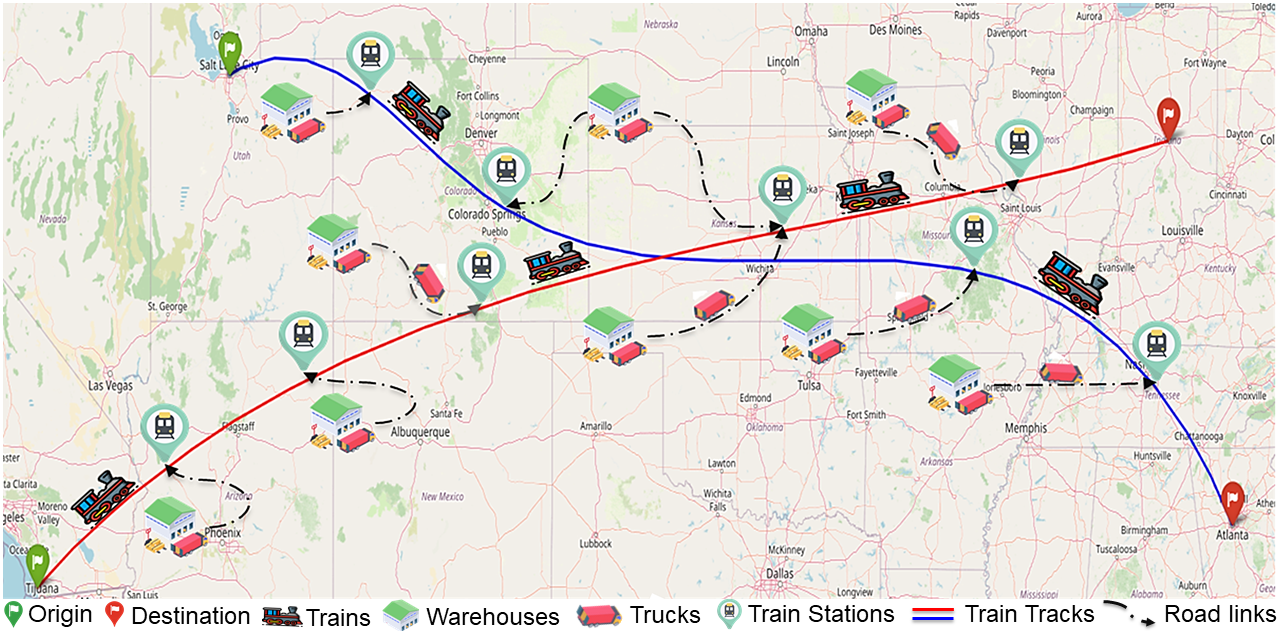}
\label{ypertaseis}
%\caption{Overview of the intermodal transport system with road-rail modal shift.}\label{fig:infographics}
\caption{Overview of the intermodal transport system: Demonstrating the strategic integration and modal shift between road and rail networks for container transport}\label{fig:infographics}
\end{figure}

This study focuses on a two-mode intermodal system (road and rail) where containers are first transported by truck to intermodal hubs, then transferred to trains operating on fixed routes and schedules. Each train serves specific destinations, and the company must align container preparation and dispatch with train departures. The optimization problem involves determining how many containers to prepare at each supply location and how to allocate them across hubs, aiming to meet destination demand on time and at minimal cost.

A core complexity lies in the stochastic nature of both demand and train capacity. Demand at destination nodes vary, and the space on each train is uncertain, often influenced by external users and unpredictable operational factors \citep{demir2016green, lu2020fuzzy, wang2021optimizing}. These uncertainties demand adaptive strategies that can dynamically respond to changing conditions while ensuring logistical efficiency and cost-effectiveness.

To illustrate this complexity, consider a real-world example. A freight train travels from Seattle, WA, to Chicago, IL, with scheduled stops at Spokane, WA; Billings, MT; Fargo, ND; and Minneapolis, MN. Multiple supply facilities located near these hubs must decide how many containers to prepare and where to route them. Unused capacity leads to inefficiencies, while over-preparation may create bottlenecks. Compounding this is the need to synchronize road transport with fixed train departures, accounting for varying travel times and uncertain spot availability. These dynamics make planning a multi-layered, high-stakes decision process.

Key constraints include: fixed train schedules, uncertain capacity at hubs, and time-sensitive delivery windows. Companies must prepare containers (regardless of load completeness) and transport them via predetermined road routes to the hubs. While routing decisions between supply points and hubs are excluded from this model, companies
consider such routes are often fixed due to operational consistency based on their frequent, standardized nature in practice \citep{basallo2021planning}. Environmental impact is another important dimension. The model seeks to reduce emissions by prioritizing nearby supply locations, minimizing road distance and aligning with policy frameworks such as the Clean Air Act \citep{epa2024cleanair} and state-level carbon reduction targets \citep{wri2024usstate}.

Overall, the problem is defined by the need to balance cost, service reliability, and sustainability under high uncertainty. We adopt the following assumptions, grounded in operational norms and common practice:
i) Containers are shipped as-is once ready, regardless of load completeness, reflecting industry practices where dispatch cannot be delayed for full consolidation;
ii) supply locations are assumed to be pre-stocked with containers, isolating logistics planning from inventory sourcing;
iii) train routes and schedules, are treated as fixed, predefined by carriers and not subject to shipper control;
iv) road travel times between supply points and hubs considers traffic conditions, suitable for strategic planning where high-frequency variations are averaged out;
v) containers are uniform in size, per ISO standards, but may vary in weight due to partial loading
vi) containers remain onboard until their final destination, avoiding mid-route handling, reflecting common linehaul practices; and
vii) road routing between supply points and hubs is predetermined, consistent with fixed or established distribution protocols \citep{basallo2021planning}

This problem setup reflects the real-world complexity of intermodal freight systems, shaped by uncertain capacities, tight scheduling, and environmental constraints. The following section presents a two-stage stochastic optimization model that captures these dynamics and supports adaptive, cost-effective planning.

\section{Optimization Model}
\label{Optimization Model}
This section presents a two-stage stochastic optimization model designed to address intermodal transportation challenges, particularly uncertainties in demand and train capacities. The model integrates a CVaR framework to balance cost efficiency and risk management. Key model components, including sets, parameters, and variables, are outlined in Table~\ref{Table:Combined}. The formulation is explained alongside a discussion of constraints, and a deterministic equivalent of the stochastic model is provided.

\begin{table}[!ht]
	\caption{Sets, Parameters, and Variables Definitions}\label{Table:Combined}
	\begin{center}
		\begin{tabularx}{\textwidth}{l X}
			\hline
			\textbf{Sets and Indices} & \textbf{Definition} \\ \hline
                $I$ $(i)$ & Origin locations (warehouses/manufacturing centers) \\
                $J$ $(j)$ & Train stations/intermodal hubs \\
                $T$ $(t)$ & Time periods \\
                $\Omega$ $(\omega)$ & Scenarios \\
                $N$ $(n)$ & Trains \\ \hline
			\textbf{Parameter} & \textbf{Definition} \\ \hline
			$C_{ij}$ & Base transportation cost per container from origin $i$ to station $j$ \\
			$\tau_{ij}$ & Travel time from origin $i$ to station $j$ \\
			$\alpha$ & Confidence level for CVaR \\
			$\lambda$ & Trade-off coefficient between cost and risk \\
			$\theta$ & Value-at-Risk (VaR) threshold for cost \\
			$\rho$ & Penalty for excess emissions \\
			$\pi_\omega$ & Unit penalty cost for missing the train schedule under scenario $\omega$ \\
			$\epsilon$ & Total allowable carbon emission for planning horizon \\
			$E_{t}$ & Carbon emission per container per time period $t$ \\
			$K_{jn\omega}$ & Random capacity of train $n$ at station $j$ in scenario $\omega$ \\
			$O_i$ & Operational cost for preparing containers at origin $i$ \\
			$p_\omega$ & Probability of scenario $\omega$ \\
			$D_{n\omega}$ & Demand for train $n$ in scenario $\omega$ \\
			$\sigma_{nj}$ & Predefined schedule, specifying the time at which train $n$ will be at station $j$ \\
			$\kappa_i$ & Maximum number of containers that can be prepared at origin $i$ \\ 
                $M$ & A large positive constant \\\hline
			\textbf{Decision Variable} & \textbf{Definition} \\ \hline
            $x_{ijn\omega t}$ & Containers shipped from origin $i$ to station $j$ for train $n$ at time $t$ in scenario $\omega$ \\
            $y_i$ & Containers prepared at origin $i$ \\
            $z_{ijn\omega t}$ & Binary variable; equals 1 if a container is assigned from $i$ to $j$ for train $n$ at time $t$ in scenario $\omega$, 0 otherwise \\ \hline
            \textbf{Auxiliary Variable} & \textbf{Definition} \\ \hline
            $U_{n\omega}$ & Unmet demand for train $n$ in scenario $\omega$ \\
            $\eta_\omega$ & Excess emissions in scenario $\omega$ \\
            $\mathcal{I}_{jn\omega}$ & Inventory of containers on train $n$ at station $j$ in scenario $\omega$ \\ \hline
		\end{tabularx}
	\end{center}
\end{table}

In our proposed model, we introduce a two-stage stochastic programming framework that integrates risk aversion, where both the objective function and second stage constraints incorporate stochastic parameters in line with acceptable practice. Specifically, we adopt a mean-risk objective that incorporates CVaR with $\alpha \in (0, 1)$ being the confidence level. This choice of risk measure allows us to model the risk preferences of decision makers concerning the stochastic cost, offering flexibility to accommodate various thresholds from no risk ($\alpha = 0$) to worst-case scenarios ($\alpha \rightarrow 1$).

To ensure the robustness of our approach, we control unfulfilled orders using demand constraints by penalizing these violations in the objective function. The first-stage problem in the model determines the operational cost incurred for preparing the containers and the number of containers shipped from each supply location based on the capacity of each facility. The problem in the second-stage focuses on meeting the demand and utilizing the available capacity on each train at each station. From this perspective, our model aims to design a reliable network for logistics operations by imposing constraints that ensure demand is fulfilled with a high specified probability.

Let the first-stage decision vector be represented by $(a, b)$, the stochastic data vector $\zeta(\omega)$, and the second-stage problem with state $\gamma(\omega)$. The optimal stochastic objective value for the second-stage be denoted as $S''(a, b, \zeta(\omega), \gamma(\omega))$.

Given the discussion above the objective of the second stage problem is stated as: 
\vspace{-0.8em}
\begin{equation}
S''(a, b, \zeta(\omega), \gamma(\omega)) = \min \left(\sum_{i \in I} \sum_{j \in J} \sum_{n \in N} \sum_{t \in T} C_{ij} x_{ijn\omega t} + \pi_\omega \sum_{n \in N} U_{n\omega}\right) \label{eq:2nd_stage_obj}
\end{equation}

The objective function in \eqref{eq:2nd_stage_obj} minimizes total transportation cost from origins to train stations and penalized unmet demand for each train in every scenario. This dual-purpose objective ensures efficient logistics operations by optimizing transportation routes and mitigating the impact of any shortfalls in meeting the demand for each train in every scenario.

Choices are typically made before uncertainties are realized, with recourse decisions deferred to the second stage that addresses the stochasticity. Once the uncertainty is known, detailed second-stage distribution decisions can be addressed. Thus, first-stage logistics planning concentrates on the number of containers prepared at each facility, while the second-stage problem manages demand and capacity uncertainties.

\subsection{Conditional Value-at-Risk (CVaR)}

Before presenting the full stochastic optimization model, we briefly review the CVaR risk measure, which plays a central role in our formulation. Originally introduced by \cite{rockafellar2000optimization} for portfolio optimization, CVaR has been widely adopted in operations research and logistics due to its ability to quantify the risk of extreme outcomes beyond the traditional Value-at-Risk (VaR) threshold. Applications of CVaR in stochastic programming include two-stage recourse formulations \citep{kunzi2006computational} and efficient polyhedral approximations for large-scale optimization problems \citep{fabian2008handling}. For better understanding see \citep{kunzi2006computational, fabian2008handling}

At a given confidence level $\alpha \in (0, 1)$, CVaR measures the expected value of losses exceeding the corresponding VaR threshold. Mathematically, CVaR can be defined as:
\vspace{-0.8em}
\begin{align}
\text{CVaR}_{\alpha}(R) = \min_{\text{VaR} \in \mathbb{R}} \left\{ \text{VaR} + \frac{1}{1 - \alpha} \mathbb{E} \left( [R - \text{VaR}]_+ \right) \right\} \label{eq:CVaR1}
\end{align}

Here, $R$ is a random variable representing the loss, and $[m]_+ = \max(m, 0)$ denotes the positive part of $m$. The VaR represents the $\alpha$-quantile of $R$, and CVaR captures the average loss beyond this quantile. Intuitively, CVaR focuses not just on the threshold where extreme losses begin, but on the full tail distribution of worse-than-expected outcomes, making it a more robust risk metric for uncertain environments.

To facilitate computational tractability, \eqref{eq:CVaR1} can be reformulated as a linear program. Introducing an auxiliary variable $\theta$ to represent the VaR estimate, and slack variables $\pi_\omega$ to model excess losses in each scenario $\omega \in \Omega$, the equivalent formulation becomes:
\vspace{-0.8em}
\begin{align}
\text{CVaR}_{\alpha}(R) = \min_{\theta, \pi} \left\{ \theta + \frac{1}{1 - \alpha} \sum_{\omega \in \Omega} p_\omega \pi_\omega \mid \pi_\omega \geq y_\omega - \theta, \; \pi_\omega \geq 0, \; \forall \omega \in \Omega \right\} \label{eq:CVaR2}
\end{align}

In this expression, $\theta$ approximates the VaR at level $\alpha$, $\pi_\omega$ captures the excess loss in scenario $\omega$ beyond $\theta$, $p_\omega$ denotes the probability of scenario $\omega$ occurring, $y_\omega$ represents the realization of $R$ under scenario $\omega$.

This linear programming form allows CVaR minimization to be embedded efficiently into stochastic optimization models, such as our two-stage framework for intermodal logistics.

Further, CVaR admits a dual representation that highlights its interpretation as a worst-case expectation over adverse outcomes. The dual form can be expressed as:
\vspace{-0.8em}
\begin{align}
\text{CVaR}_\alpha(R) = \max_{q} \left\{ \frac{1}{1 - \alpha} \sum_{\omega \in \Omega} q_\omega y_\omega \mid \sum_{\omega \in \Omega} q_\omega = 1 - \alpha, \; 0 \leq q_\omega \leq p_\omega, \; \forall \omega \in \Omega \right\} \label{eq:CVaR3}
\end{align}

Alternatively, it can be interpreted through an integral form over the upper tail quantiles:
\vspace{-0.8em}
\begin{align}
\text{CVaR}_\alpha(R) = \frac{1}{1 - \alpha} \int_{\alpha}^1 \theta_q(R) \, dq
\end{align}

where $\theta_q(R)$ denotes the $q$-quantile of the loss distribution.

The dual form emphasizes that CVaR aggregates weighted worst-case losses beyond the VaR threshold, providing a coherent and convex risk measure, satisfing desirable properties such as monotonicity, subadditivity, and translation equivariance \citep{artzner1999coherent}. This focus on extreme scenarios makes CVaR well-suited for applications like intermodal logistics, where demand variability and capacity fluctuations can lead to significant operational disruptions.

This dual structure enables flexible integration into optimization models, supporting resilient logistics planning under extreme uncertainties. In our work, CVaR is incorporated directly into the objective function of the two-stage stochastic optimization model, allowing for simultaneous minimization of expected transportation costs and explicit control over downside risk. This risk-aware formulation enhances the model's ability to generate resilient, cost-effective, and environmentally sustainable operational plans. Building on this foundation, we now develop the two-stage stochastic optimization model that incorporates CVaR for intermodal freight logistics under uncertainty.

\subsection{The intermodal stochastic optimization model}

To address the complexities of managing intermodal logistics under uncertainty, we present a proposed model that integrates cost efficiency, risk mitigation, and environmental sustainability. This formulation emphasizes a dual-stage decision-making process, where initial resource allocation is followed by adaptive adjustments based on real-time conditions. By structuring the model in this way, it captures the dynamic interplay between strategic planning and operational flexibility, setting the foundation for the detailed objectives and constraints discussed below.
\subsubsection{Objective}
\vspace{-1.2em}
\begin{align}
\text{Min} & \left( \sum_{i \in I} O_i y_i + (1 - \lambda) \mathbb{E} \left(S''(a, b, \zeta(\omega), \gamma(\omega))) \right) + \lambda \text{CVaR}_{\alpha} \left( S''(a, b, \zeta(\omega), \gamma(\omega)) \right) + \rho \sum_{\omega \in \Omega} p_\omega \eta_\omega \right)\label{eq:full_obj}
\end{align}

The objective function \eqref{eq:full_obj} integrates both the preparation (first-stage) and operational costs (second-stage)  for the intermodal stochastic optimization model. This objective minimizes the mean-risk function related to the total stochastic cost. Here, ($\sum_{i \in I} O_i y_i$) represents the fixed cost of preparing containers at origin locations. The term $(1 - \lambda) \mathbb{E} \left( S''(a, b, \zeta(\omega), \gamma(\omega)) \right)$ reflects the expected second-stage cost, including transportation and unmet demand across all scenarios. The term $\lambda \text{CVaR}_{\alpha} \left( S''(a, b, \zeta(\omega), \gamma(\omega))) \right)$ captures tail-risk aversion, where $\lambda$ and $\alpha$ are risk parameters controlling the weight and confidence level, respectively. The final term, $\rho \sum_{\omega \in \Omega} p_\omega \eta_\omega$, penalizes emissions exceeding allowable thresholds, thereby embedding sustainability directly into the decision-making process.
 
\subsubsection{Constraints}
\vspace{-0.8em}
\begin{align}
    &\sum_{j \in J} \sum_{n \in N} \sum_{t \in T} x_{ijn\omega t} \leq y_i, &&\forall i \in I, \, \forall \omega \in \Omega \label{eq:supply}\\
    &\sum_{i \in I} \sum_{t \in T} x_{ijn\omega t} \leq K_{jn\omega}, &&\forall j \in J, \, n \in N, \, \omega \in \Omega \label{eq:train_capacity}\\ 
    &x_{ijn\omega t} \leq M z_{ijn\omega t}, &&\forall i \in I, \, j \in J, \, n \in N, \, t \in T, \, \omega \in \Omega \label{constraint:linking}\\
    &\sum_{i \in I} \sum_{j \in J} \sum_{n \in N} \sum_{t \in T} x_{ijn\omega t} E_{t} \tau_{ij}\leq \epsilon + \eta_\omega, &&\forall \omega \in \Omega \label{eq:emissions}\\
    &\mathcal{I}_{j_1n\omega} = \sum_{i \in I} \sum_{t \in T} x_{ij_1n\omega t}, &&\forall n \in N, \, \omega \in \Omega \label{eq:initial_inventory}\\
    &\mathcal{I}_{jn\omega} = \mathcal{I}_{(j-1)n\omega} + \sum_{i \in I} \sum_{t \in T} x_{ijn\omega t}, &&\forall j \in J \setminus \{j_1\}, \, n \in N, \, \omega \in \Omega \label{eq:cumulative_inventory}\\
    & \sum_{i \in I} \sum_{j \in J} \sum_{t \in T} x_{ijn \omega t} + U_{n\omega} = D_{n\omega}, && \forall n \in N, \omega \in \Omega \label {eq:demand}\\
    &\mathcal{I}_{j_fn\omega} \leq D_{n\omega}, &&\forall n \in N, \, \forall \omega \in \Omega \label{eq:final_inventory}\\
    & \sum_{i \in I} x_{ijn \omega t} \leq y_i, && \text{if }t + \tau_{ij} \leq \sigma_{nj} \quad \forall j \in J, n \in N, \omega \in \Omega, t \in T \label{eq:travel_time}\\ 
    & \sum_{i \in I} x_{ijn \omega t} = 0 && \text{if }t + \tau_{ij} > \sigma_{nj} \quad \forall j \in J, n \in N, \omega \in \Omega, t \in T \label{eq:combined_zero_flow} \\ %combines the two above
    & y_i \leq \kappa_i, &&\forall i \in I \label{eq:upper_bound_relationship} \\
    &z_{ijn\omega t} \in \{0, 1\}, &&\forall i \in I, \, j \in J, \, n \in N, \, t \in T, \, \omega \in \omega \label{eq:binary_penalty}\\
    &x_{ijn\omega t},y_i,\mathcal{I}_{jn\omega}, U_{n\omega}  \in \mathbb{Z}_{\geq 0}, \quad \sigma_{nj}, \eta_\omega \geq 0, &&\forall i \in I, \, j \in J, \, n \in N, \, t \in T, \, \omega \in \Omega \label{eq:non_negativity}
\end{align}

The constraints applied to the model effectively determine demand assignments, the number of containers on each train per hub, and limit the carbon emissions during the time period.
To start with, constraints \eqref{eq:supply} ensure that the total number of containers transported from each origin $i$ does not exceed the available supply at the origin.  Constraints \eqref{eq:train_capacity} ensure that the number of containers transported to each train station does not exceed the available capacity on each train at a particular station. Further, constraints \eqref{constraint:linking} ensure that container flows \( x_{ijn\omega t} \) are linked to  binary decision variables \( z_{ijn\omega t} \), ensuring that containers can only be shipped if the corresponding binary variable is activated (\( z_{ijn\omega t} = 1 \)). The large constant \( M \) acts as an upper bound, preventing any flow \( x_{ijn\omega t} \) unless the route is selected.

As for constraints \eqref{eq:emissions}, they ensure that the total emissions do not exceed the specified emissions cap plus any allowable excess. Constraints \eqref{eq:initial_inventory} initialize the inventory at the first station for each train and scenario, meaning the number of containers on the train at the first station will be the number of containers that were transported to that station. While constraints \eqref{eq:cumulative_inventory} update the inventory at each subsequent station based on the previous station's inventory and the incoming flow of containers, constraints \eqref{eq:demand} ensure that the total number of containers transported to each train station and the unmet demand sum up to the total demand for each train. Constraints \eqref{eq:final_inventory} ensure that the inventory at the final station \( j_f \) for each train across all scenarios does not exceed the demand for that train.  Constraints \eqref{eq:travel_time} dictate that containers from origin  are shipped to a train station at a certain time only if the shipping time plus the travel time  does not exceed the train's scheduled departure time.

We make sure that no containers are transported from any origin to a station at a particular time if the container cannot meetup with the train schedule (see Constraints \eqref{eq:combined_zero_flow}). Constraints \eqref{eq:upper_bound_relationship} ensure that the available supply at each origin $i$, does not exceed the maximum allowable supply at each origin. Lastly, the binary penalty constraints \eqref{eq:binary_penalty} and the non-negativity constraints \eqref{eq:non_negativity} ensure logical consistency and enforce the non-negativity of the decision variables.

\subsection{Deterministic equivalent of stochastic optimization model}
To simplify implementation and analysis, the stochastic model is reformulated into its deterministic equivalent. This reformulation explicitly separates the CVaR term, improving the interpretability of risk-related costs and facilitating the identification of their individual contributions to the overall objective. Relating our stochastic formulation to the knapsack-style approach outlined in \eqref{eq:CVaR3}, we define the random total cost $R$ as:
\vspace{-1.5em}

\begin{equation}
    R = \left( \sum_{i \in I} O_i y_i + S''(a, b, \zeta(\omega), \gamma(\omega)) + \rho \eta_\omega \right)
\end{equation}
where $S''(a, b, \zeta(\omega), \gamma(\omega))$  represents the optimal value of the second-stage problem for scenario $\omega$, and $\rho \eta_\omega$ captures the penalties for exceeding emissions limits.

In this context, the Value-at-Risk at confidence level $\alpha$ denoted $\text{VaR}_{\alpha}(R)$, establishes an upper threshold on the total cost that is exceeded with a probability of at most $1 - \alpha$. However, $\text{CVaR}_{\alpha}(R)$ provides deeper insight by quantifying the expected total cost in scenarios where this threshold is surpassed. This measure captures the severity of extreme cost events, emphasizing potential losses beyond the VaR threshold.

Building on this risk assessment framework, we now present the reformulated deterministic equivalent of our stochastic model. The objective function in \eqref{eq:full_obj} becomes:
\vspace{-0.8em}
\begin{align}
\min & \sum_{i \in I} O_i y_i + (1 - \lambda) \sum_{\omega \in \Omega} p_\omega \left( \sum_{i \in I} \sum_{j \in J} \sum_{n \in N} \sum_{t \in T} C_{ij} x_{ijn\omega t} + \pi_\omega \sum_{n \in N} U_{n\omega} \right) \nonumber \\
& + \lambda \text{CVaR}_\alpha \left( \sum_{i \in I} \sum_{j \in J} \sum_{n \in N} \sum_{t \in T} C_{ij} x_{ijn\omega t} + \pi_\omega \sum_{n \in N} U_{n\omega} \right) + \rho \sum_{\omega \in \Omega} p_\omega \eta_\omega \label{eq:objective}
\end{align}

Subject to:
\begin{flalign}
\quad \quad \quad (6) - (18) && \label{eq:con_ref}
\end{flalign}

In addition to the reformulation above with objective \eqref{eq:objective}, which explicitly includes CVaR, the model can be equivalently reformulated for clarity and highlighting separation of terms. This alternative formulation separates the CVaR component from the main summation of costs, providing a more interpretable representation of the objective. By isolating the CVaR term, it enhances the understanding of its individual contribution to the total cost and its role in risk management. This alternative form is particularly useful for both theoretical analysis and practical implementation, as it allows for a clearer distinction between the expected costs and the risk-related costs. To facilitate this understanding, we present the alternative objective function as follows:
Alternatively, \eqref{eq:objective} can be written as:
\vspace{-0.8em}
\begin{align}
\min & \sum_{i \in I} O_i y_i + (1 - \lambda) \sum_{\omega \in \Omega} p_\omega \left( \sum_{i \in I} \sum_{j \in J} \sum_{n \in N} \sum_{t \in T} C_{ij} x_{ijn\omega t} + \pi_\omega \sum_{n \in N} U_{n\omega} \right) \nonumber \\
& + \lambda \text{CVaR} + \rho \sum_{\omega \in \Omega} p_\omega \eta_\omega \label{eq:objective_CVar}
\end{align}

Also subject to (6) - (18) and the constraints below:
\begin{align}
    & \xi_\omega \geq c_{e \omega} - \theta, \quad &&\forall \omega \in \Omega \label{eq:CVaR_E1} \\
    & \xi_\omega \geq 0, \quad &&\forall \omega \in \Omega \label{eq:CVaR_E2}\\
    & \text{CVaR} = \theta + \frac{1}{1 - \alpha} \sum_{\omega \in \Omega} p_\omega \xi_\omega \label{eq:CVaR_E3}
\end{align}

Where $\xi_\omega$ denote the excess cost by which the cost exceeds the $\theta$ and $c_{e \omega}$ represent the total transportation and penalty costs for every scenario. \eqref{eq:CVaR_E1} ensures that the excess cost is calculated correctly based on the difference between the total cost in each scenario and the VaR. Constraint \eqref{eq:CVaR_E2} ensures that excess is non-negative. While \eqref{eq:CVaR_E2} defines CVaR as the VaR ($\theta$) plus the expected excess cost over all scenarios, scaled by $\frac{1}{1 - \alpha}$. CVaR represents the expected value of the excess cost beyond the $\theta$, providing a measure of the average loss in the worst $(1 - \alpha)\%$ of cases.

It is important to note that lower values of risk measures are preferred, aligning with the objective of minimizing potential extreme costs. The parameter $\lambda \in \{0, 1\}$ serves as a weighting factor that balances the trade-off between the expected cost and the risk associated with cost variability. A higher $\lambda$ 
 places more emphasis on risk aversion, potentially at the expense of higher expected costs.

In summary, the deterministic equivalent formulation provides a practical and computationally efficient framework that captures both the expected operational costs and the risks associated with uncertainties in demand and capacity. By explicitly incorporating Conditional Value-at-Risk (CVaR) into the objective function and constraints, the model empowers decision-makers to balance cost efficiency with risk mitigation, enabling the development of robust and resilient intermodal transportation plans. The next section illustrates the model's applicability using real-world data from U.S. freight networks, highlighting its ability to address operational complexities effectively and offering practical insights through a detailed case study and computational experiments.

\section{Case Study and Computational Results} \label{Computational Experiments}

\begin{figure}[!htp]
    \centering
    \captionsetup{justification=centering} % Center captions and subcaptions
    \begin{subfigure}[b]{0.48\textwidth}
        \centering
        \includegraphics[width=\textwidth, height=5cm]{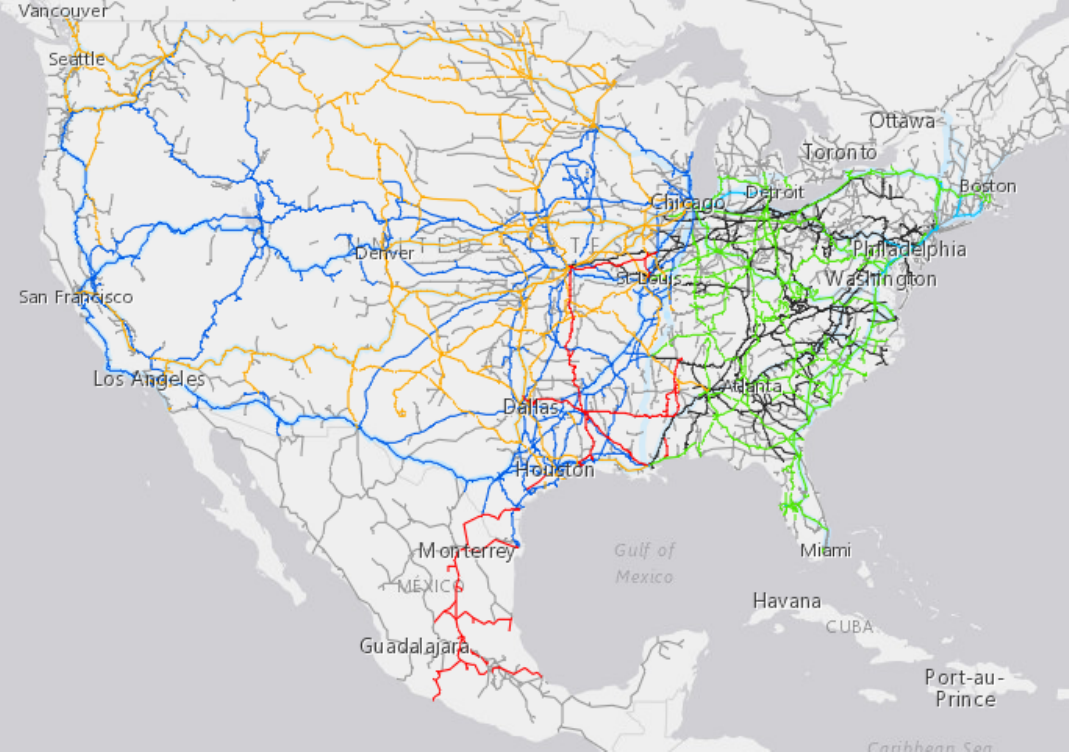}
        \caption{Overview of the rail tracks connections in the U.S.}
        \label{fig:rail_map}
    \end{subfigure}
    \hfill
    \begin{subfigure}[b]{0.48\textwidth}
        \centering
        \includegraphics[width=\textwidth, height=5cm]{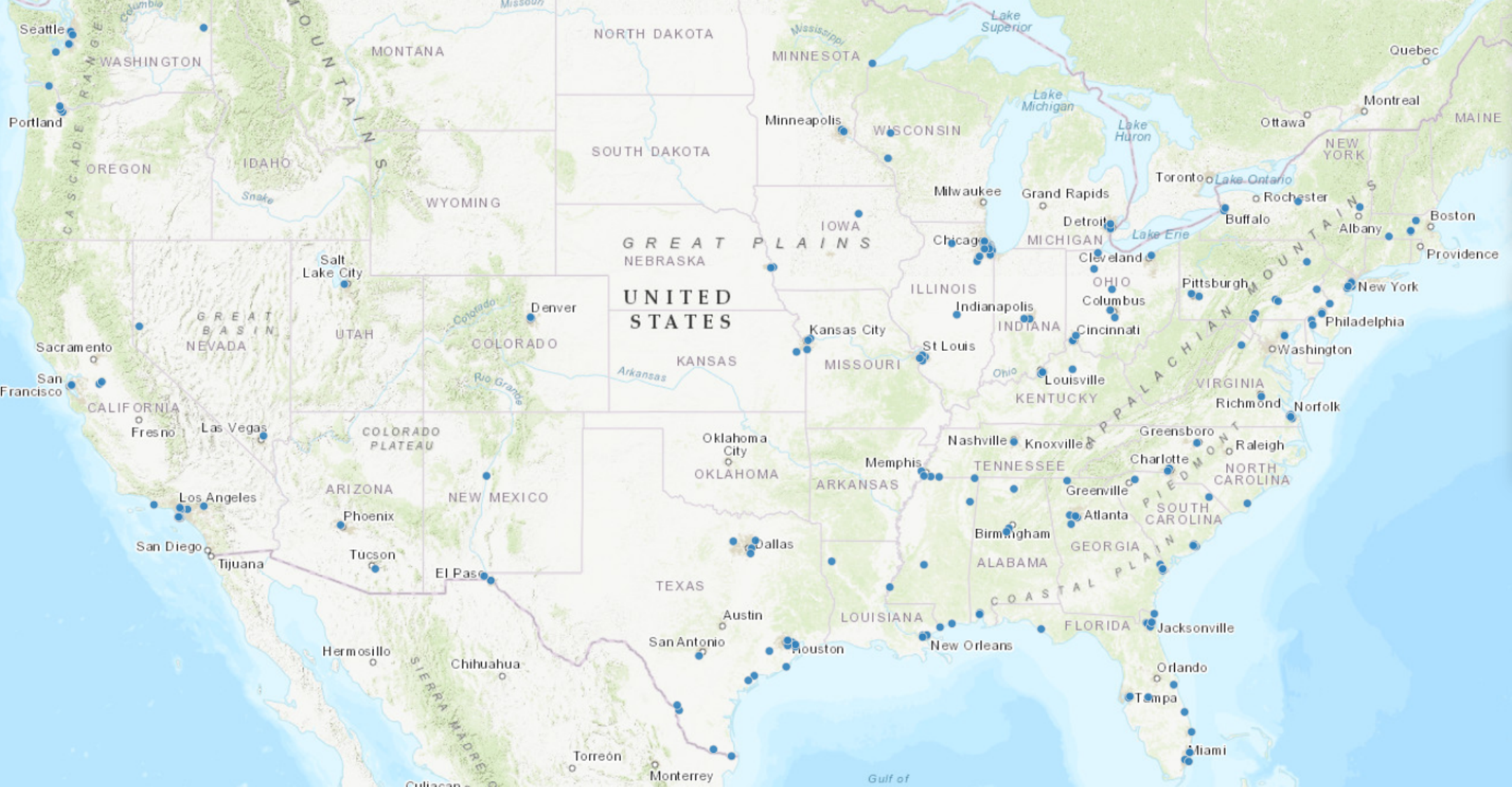}
        \caption{Intermodal hubs for freight rail in the U.S.}
        \label{fig:hubs_map}
    \end{subfigure}
    \caption{Rail infrastructure and intermodal hubs in the United States}
    \label{fig:main_figure}
\end{figure}

This section applies the proposed two-stage stochastic optimization model to an intermodal transportation problem in the U.S. The analysis utilizes publicly available datasets, including FAF \cite{bts_faf} and NTAD \cite{bts_ntad}, which provide detailed information on rail connections, intermodal terminals, and road links to hubs from supply locations such as warehouses and manufacturing facilities. To address data gaps, missing supply-hub connections were supplemented by manually retrieving data from Google Maps. Figures \ref{fig:rail_map} and \ref{fig:hubs_map} provide a visual overview of the spatial layout of the U.S. rail network, highlighting connectivity between supply points, intermodal hubs, and key corridors. This foundational data supports the development of a network model, detailed subsequently.

The computational experiments were executed on a machine equipped with a 12th Gen Intel(R) Core(TM) i7-1265U processor, 16 GB RAM, and a 64-bit operating system. Each instance was solved independently using the Gurobi Optimizer 11.0.0, ensuring efficient handling of the model's complexity. The implementation was performed in {Python 3.11.4}. The results of these experiments are detailed in the subsequent subsections.

\subsection{Policy Questions}
To guide the analysis and validate the applicability of the proposed model, the case study focuses on addressing the following research questions. These questions are designed to evaluate the model's robustness, scalability, and practicality in handling real-world intermodal transportation challenges:

\begin{enumerate}
    \item \textbf{Balancing Cost, Risk, and Decision-Making Under Uncertainty:}  
    How does the model manage trade-offs between cost efficiency and risk mitigation under varying levels of risk aversion and confidence?
    What are the implications for system resilience, cost structure, and stakeholder priorities?
    How does increasing risk aversion affect the stability and predictability of outcomes across scenarios?

    \item \textbf{Trade-Offs Between Cost and Sustainability:}  
    What are the trade-offs between cost efficiency and sustainability goals, particularly under varying levels of emissions constraints?

    \item \textbf{Adaptability to Uncertainty:}  
    How robustly does the model adapt to uncertainty in demand, capacity, and operational disruptions such as delays or unforeseen changes?

    \item \textbf{Scalability Across Networks:}  
    Can the model efficiently accommodate larger intermodal networks, including scenarios with additional supply nodes, hubs, and complex routing requirements?

    \item \textbf{Impact of Intermodal Transfers:}  
    To what extent do intermodal transfer times and associated costs influence key optimization metrics such as total cost, scheduling efficiency, and risk measures?

    \item \textbf{Impact of Capacity Constraints on Cost and Demand Satisfaction:}  
    How do train capacity affect total system costs and unmet demand, and at what point do additional capacity increases yield diminishing returns?

    \item \textbf{Impact of Risk Preferences on Cost Stability and Scenario Variability:}  
    How does increasing $\lambda$ affect cost outcomes across demand scenarios? 
    Do additional risk-mitigation efforts yield diminishing returns?
\end{enumerate}

By addressing these questions, the case study evaluates the model's ability to optimize intermodal logistics under uncertainty while balancing cost minimization, risk resilience, and sustainability objectives. The subsequent sections provide detailed analyses into these aspects, demonstrating the model’s practical applicability and contributions.

\subsection{Case Study Setup}

In this section, we detail the data used for our case study, which focuses on major freight train routes that connect 26 intermodal hubs with 50 supply locations, primarily warehouses. These routes were carefully selected to reflect the complexity of intermodal transportation in the U.S., capturing variation in demand, train capacity, and schedules. The train routes and intermodal hubs analyzed are predominantly operated by the Norfolk Southern Corporation, a major freight railroad company in the U.S. However, five of the hubs are jointly used by Norfolk Southern and Union Pacific Railroad, reflecting real-world infrastructure sharing. These shared hubs introduce coordination challenges such as scheduling and capacity constraints, which are addressed in the optimization model to enhance its applicability.

To ensure uniformity across the rail network, FedEx warehouses were used as supply locations. This choice was made to standardize the analysis and maintain consistency in evaluating logistics across different regions. FedEx warehouses, known for their strategic distribution capabilities, provide a practical and representative basis for modeling supply points in intermodal transportation. Furthermore, some of these warehouses serviced multiple intermodal hubs, reflecting the flexibility and extensive reach of distribution networks in real-world logistics.

Parameter values are drawn from peer-reviewed studies and public sources to ensure realism and relevance. $C_{ij}$ and $\pi_\omega$ follow intermodal cost structures in \citep{goodarzi2024evaluating}, while $O_i$ is set to \$250/container, similar to estimates in \citep{fhwa2025map21}. The emissions cap $\epsilon$ adopts the approach in \citep{heinold2020emission}, representing the share of emissions a shipper is willing to tolerate, and is varied in our analysis. $E_t$ uses modal emission factors provided by EPA guidance \citep{epa2024cleanair}. The emissions penalty $\rho$ is set to \$20 per metric ton of CO\textsubscript{2}, reflecting typical social cost benchmarks. Additionally, to enhance realism, the planning horizon (60 hours) is aligned with \citep{ke2021framework}, and a 1-hour per-container handling time at hubs is based on estimates from \citep{assadipour2016toll}.

\subsubsection{Network Dynamics}

This study models an intermodal freight transport network comprising five scheduled train routes spanning major U.S. economic corridors. Each route connects strategic hubs via fixed schedules, facilitating container transfers between rail and road. Train operations are characterized by defined origin–destination pairs, intermediate stops, and fixed dwell times of approximately six hours per stop to allow for loading, unloading, and intermodal handling.

As an illustrative example, Train A travels from Buffalo, NY, to Detroit, MI (Delray), covering ~2,105 miles. It stops at Greencastle, Columbus, Cincinnati, St. Louis, and Decatur. Each segment has calibrated travel times based on actual rail distances, with intermediate dwell periods to align with realistic freight operations. The routing structure reflects real-world transit logistics and timing constraints. For instance, the Buffalo–Greencastle segment spans 280 miles and 6.5 hours, while Decatur–Detroit covers 310 miles over 7 hours.

Each hub is linked to multiple nearby facilities (e.g., warehouses or plants), with distances typically under 100 miles. These are selected based on proximity and relevance to freight demand. For example, FedEx Warehouse 13 is 68 miles from Buffalo, while facilities near St. Louis range from 11 to 29 miles away. These first-mile connections are assumed to use fixed road routes, aligning with typical industry practices. Road segments are modeled using average speeds (55 mph) and costs (\$3/mile), enabling precise estimation of total delivery time and cost from origin to hub.

Other routes (Trains B–E) span cities such as Seattle, Dallas, Jacksonville, and New Orleans, incorporating 4–5 intermediate hubs each, which can be recreated. While Train A is detailed for clarity, the others follow similar structures, and their inclusion demonstrates the model’s scalability across network configurations.

\subsection{Computational Experiment}

This section presents a comprehensive analysis of the proposed model's performance under various configurations and parameter settings addressing the proposed research questions. Through a series of computational experiments, we evaluate the sensitivity of key model outputs, such as total cost and risk measures, to changes in $\lambda$ (risk aversion weight) and $\alpha$ (risk confidence level). The experiments also examine the model’s ability to balance cost minimization with risk mitigation and assess its adaptability to different logistical scenarios.

To ensure clarity, several abbreviations are used throughout this section. OBJ denotes the stochastic model’s objective function, incorporating both first-stage and second-stage costs. E(TC) represents the expected total cost across all scenarios, while ASC refers to the average stochastic cost, capturing the mean realized costs across experiments. CVaR quantifies the expected cost of extreme scenarios beyond a specified threshold, determined by $\alpha$.

The parameter $\lambda$ controls trade-off between expected cost minimization and risk aversion. It is varied from 0 to 1 in increments of 0.1 (i.e., $\lambda = 0, 0.1, 0.2, \ldots, 1$), allowing us to explore the full spectrum of decision-making priorities, from cost-focused strategies ($\lambda = 0$) to highly risk-averse approaches ($\lambda = 1$). Similarly, the risk confidence level $\alpha$ is tested at values of 0.25, 0.5, 0.75, and 1. These values provide insights into the model’s response to varying levels of risk tolerance, from focusing on extreme losses in a quarter of scenarios ($\alpha = 0.25$) to worst-case outcomes ($\alpha = 1$).

The demand scenarios analyzed in this study include extreme scenarios, where demand is significantly higher or lower than other scenarios, representing fluctuations that arises due to seasonal shifts, economic disruptions, or supply chain imbalances. These variations enable an assessment of system performance under both exceptionally high and low freight volumes, capturing intermodal hubs and train stations with the highest stochasticity and examining their impact on capacity uncertainty at each station. By incorporating these extreme cases, the model captures a realistic spectrum of uncertainty, ensuring that optimization strategies remain effective across a wide range of operating conditions.

By analyzing these scenarios, we aim to demonstrate the robustness and practical applicability of the model in optimizing intermodal logistics under uncertainty. The subsequent subsections provide a detailed breakdown of the experiments, starting with an analysis of the impact of $\alpha$ and $\lambda$ on the expectation and CVaR of the total cost.

\subsubsection{Impact of Risk-Aversion Parameter (\texorpdfstring{$\lambda$}{λ}) and Confidence Level (\texorpdfstring{$\alpha$}{α}) on Cost and CVaR}

The interplay between $\lambda$ and $\alpha$ provides valuable insights into how the model balances total cost and risk. Understanding the effect of these parameters is crucial in decision-making under uncertainty. The parameter $\lambda$ dictates the weight placed on the CVaR, controlling the trade-off between minimizing average costs and reducing exposure to extreme cost scenarios. Meanwhile, $\alpha$ determines the quantile of the cost distribution used for CVaR calculations, thereby capturing the confidence level for considering risk reflecting the level of conservatism in risk management. 

\begin{figure}[!htp]
    \centering
    \captionsetup{justification=centering} % Center captions and subcaptions
    \begin{subfigure}[b]{0.48\textwidth}
        \centering
        \includegraphics[width=\textwidth, height=5cm]{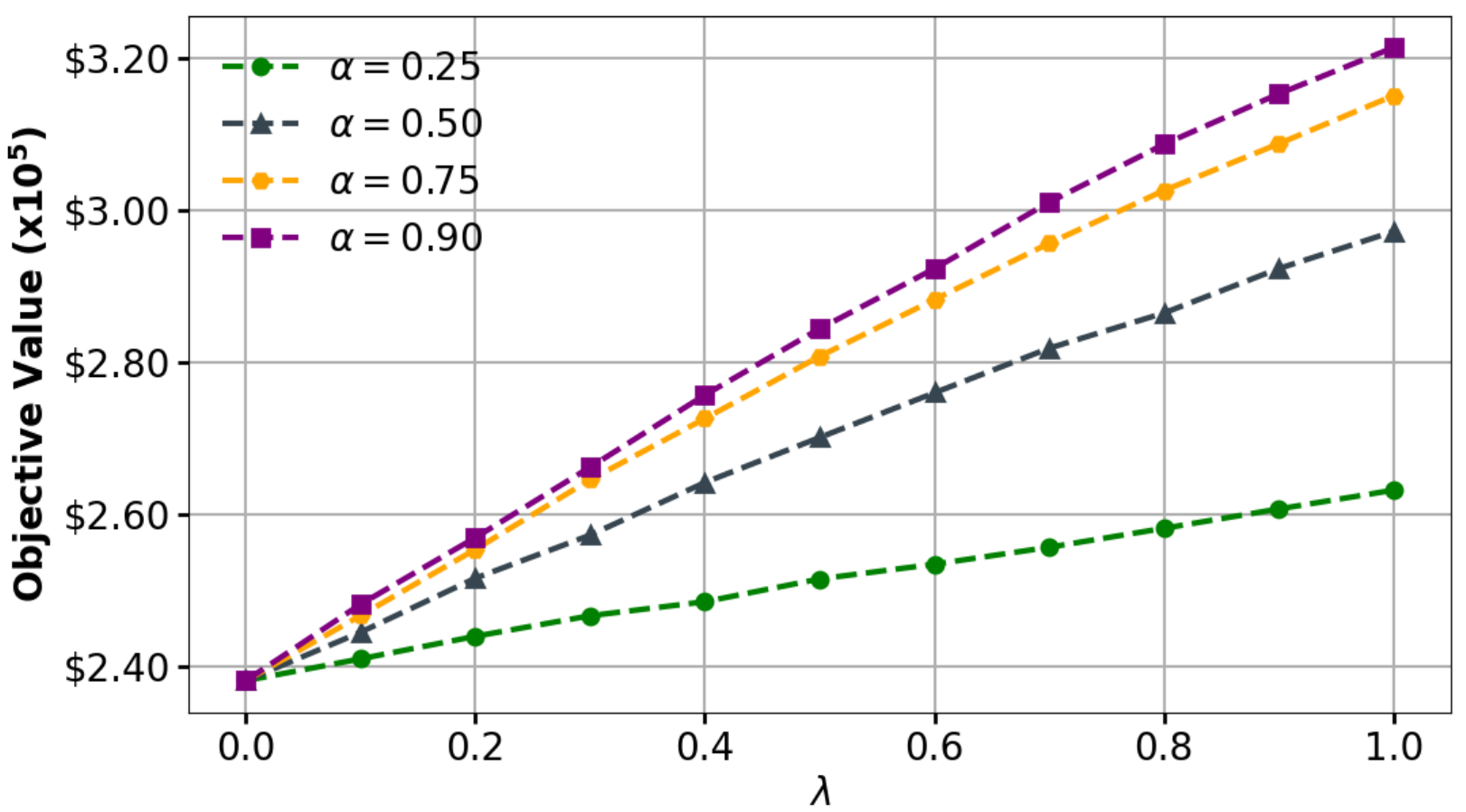}
        \caption{OBJ vs. $\lambda$ for different $\alpha$}
        \label{fig:OBJ_lambda}
    \end{subfigure}
    \hfill
    \begin{subfigure}[b]{0.48\textwidth}
        \centering
        \includegraphics[width=\textwidth, height=5cm]{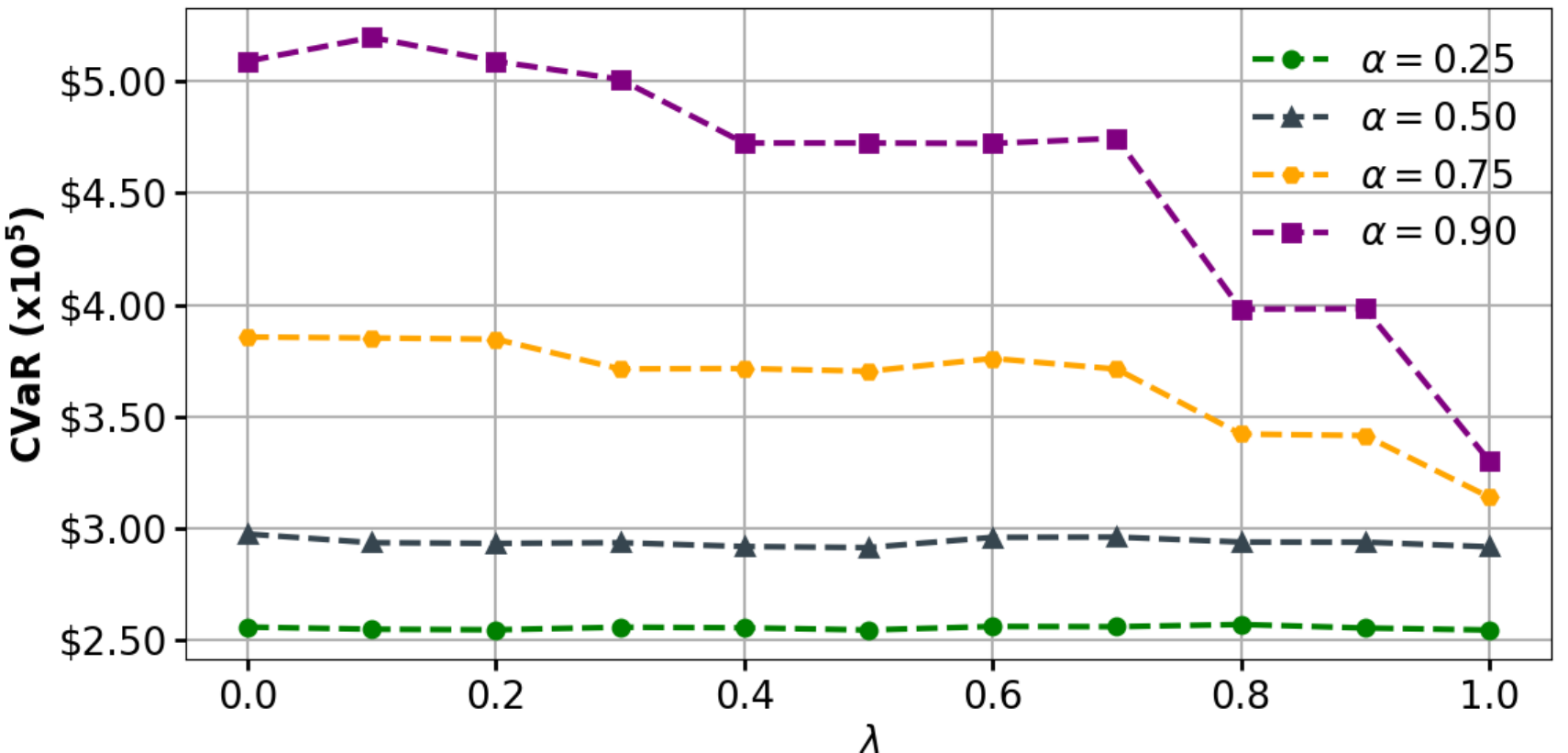}
        \caption{$\mathbf{CVaR}$ vs. $\lambda$ for different $\alpha$}
        \label{fig:CVAR_lambda}
    \end{subfigure}
    \caption{Relationship between $\lambda$, $\mathbf{OBJ}$, and $\mathbf{CVaR}$ for varying $\alpha$ levels.}
    \label{fig:OBJ_CVAR_analysis}
\end{figure}

\begin{figure}[!htp]
    \centering
    \captionsetup{justification=centering} % Center caption alignment
    \includegraphics[width=0.7\textwidth, height=6cm]{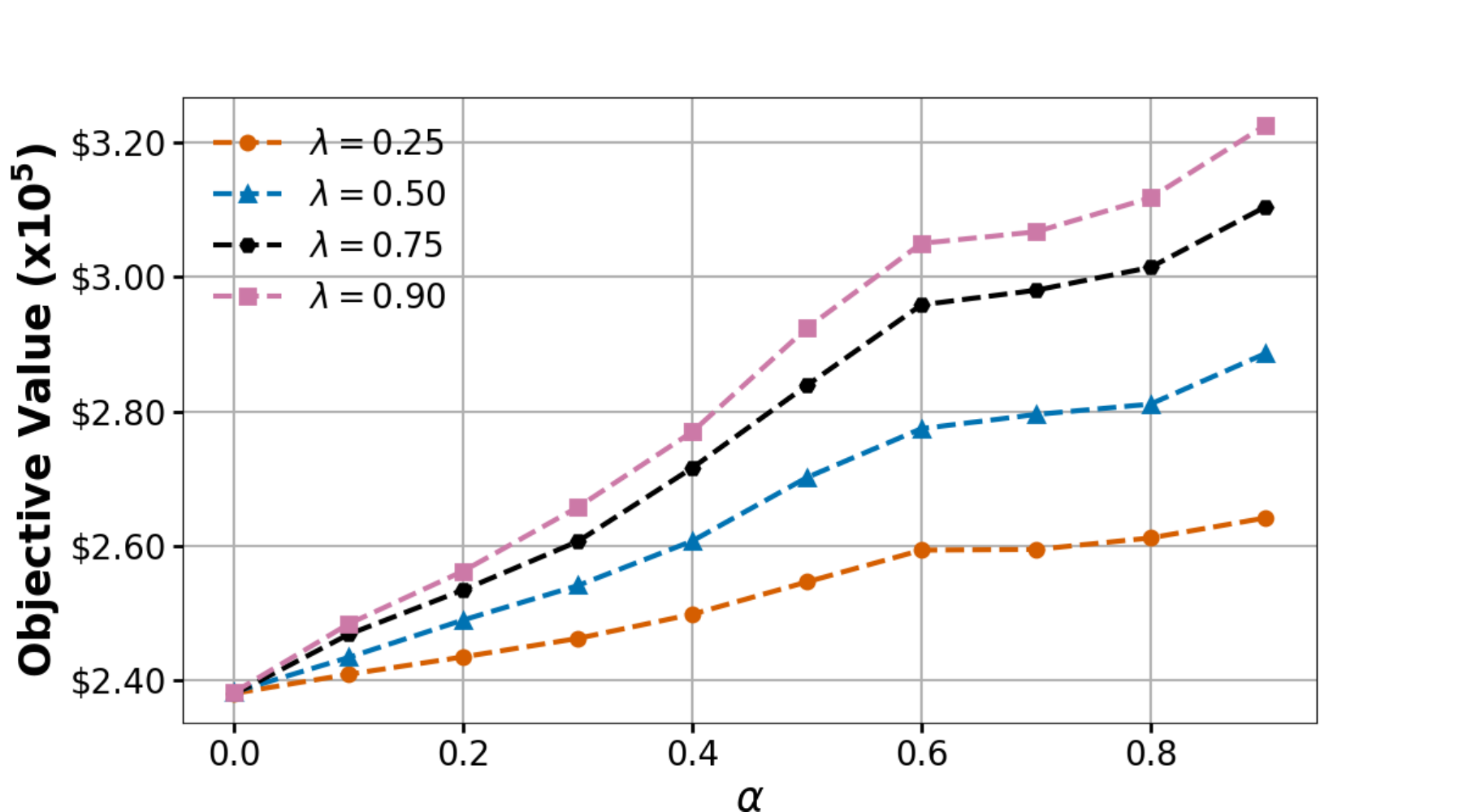}
    \caption{Impact of $\alpha$ on the Objective Value for fixed $\lambda$ values.}
    \label{fig:OBJ_vs_alpha_fixed_lambda}
\end{figure}

The combined effect of these parameters is presented in Table~\ref{tab:alpha_lambda}, summarizing the relationship between $\lambda$, $\alpha$, and key cost metrics, and further illustrated in Figures~\ref{fig:OBJ_lambda}, \ref{fig:CVAR_lambda}, and \ref{fig:OBJ_vs_alpha_fixed_lambda}. This analysis directly addresses the first policy question on balancing cost and risk by demonstrating how changes in $\lambda$ and $\alpha$ influence total cost and extreme-event mitigation. Additionally, it provides insights into how risk-based decision-making affects cost structures and system resilience.

\begin{table}[h!]
\centering
\scriptsize % Reduce font size of the table
\caption{Results for different $\alpha$ and $\lambda$ values (costs in \$ $\times$ 10\textsuperscript{3})}
\label{tab:alpha_lambda}
%\resizebox{\textwidth}{!}{%
\begin{tabular}{c|ccccc|ccccc}
\toprule
\multicolumn{1}{c}{} & \multicolumn{5}{c|}{$\alpha = 0.25$} & \multicolumn{5}{c}{$\alpha = 0.5$} \\
\cmidrule(r){1-6} \cmidrule(l){7-11}
\textbf{Lambda} & \textbf{OBJ} & \textbf{ASC} & \textbf{E(TC)} & \textbf{VaR} & \textbf{CVaR} & \textbf{OBJ} & \textbf{ASC} & \textbf{E(TC)} & \textbf{VaR} & \textbf{CVaR} \\
\midrule
0.0 & 238.2 & 210.1 & 226.5 & 150.0 & 256.0 & 238.2 & 210.1 & 226.5 & 167.7 & 297.6 \\
0.1 & 241.1 & 209.3 & 225.7 & 149.8 & 255.0 & 244.6 & 211.1 & 227.3 & 176.2 & 293.7 \\
0.2 & 244.0 & 210.9 & 225.9 & 150.6 & 254.8 & 251.7 & 214.5 & 229.2 & 176.2 & 293.4 \\
0.3 & 246.7 & 210.5 & 226.4 & 150.2 & 255.9 & 257.3 & 212.9 & 228.3 & 176.2 & 293.7 \\
0.4 & 248.6 & 214.1 & 227.4 & 150.2 & 255.7 & 264.2 & 212.1 & 226.9 & 174.2 & 292.0 \\
0.5 & 251.6 & 210.7 & 225.8 & 150.6 & 254.7 & 270.2 & 218.7 & 229.8 & 177.6 & 291.5 \\
0.6 & 253.5 & 216.7 & 229.8 & 160.2 & 256.2 & 276.0 & 220.3 & 232.7 & 178.6 & 296.1 \\
0.7 & 255.7 & 218.2 & 230.2 & 162.4 & 256.1 & 281.9 & 220.5 & 232.9 & 178.3 & 296.2 \\
0.8 & 258.2 & 216.0 & 229.6 & 155.6 & 257.2 & 286.5 & 225.8 & 236.1 & 189.5 & 294.0 \\
0.9 & 260.8 & 218.2 & 230.1 & 164.2 & 255.6 & 292.4 & 223.2 & 234.7 & 185.7 & 294.0 \\
1.0 & 263.2 & 226.6 & 234.0 & 173.8 & 254.7 & 297.2 & 240.0 & 245.1 & 202.8 & 291.9 \\
\midrule
\multicolumn{1}{c}{} & \multicolumn{5}{c|}{$\alpha = 0.75$} & \multicolumn{5}{c}{$\alpha = 0.9$} \\
\cmidrule(r){1-6} \cmidrule(l){7-11}
\textbf{Lambda} & \textbf{OBJ} & \textbf{ASC} & \textbf{E(TC)} & \textbf{VaR} & \textbf{CVaR} & \textbf{OBJ} & \textbf{ASC} & \textbf{E(TC)} & \textbf{VaR} & \textbf{CVaR} \\
\midrule
0.0 & 238.2 & 210.1 & 226.5 & 200.7 & 385.7 & 238.2 & 210.1 & 226.5 & 200.7 & 509.0 \\
0.1 & 246.8 & 212.7 & 228.3 & 200.1 & 385.2 & 248.2 & 211.1 & 227.1 & 194.4 & 519.6 \\
0.2 & 255.5 & 216.1 & 229.6 & 200.9 & 384.7 & 257.0 & 215.7 & 229.6 & 200.3 & 508.9 \\
0.3 & 264.6 & 211.6 & 228.3 & 215.1 & 371.4 & 266.3 & 211.0 & 227.7 & 201.4 & 500.9 \\
0.4 & 272.7 & 221.7 & 232.3 & 216.1 & 371.6 & 275.8 & 216.4 & 230.9 & 216.1 & 472.4 \\
0.5 & 280.8 & 219.4 & 231.8 & 216.1 & 370.4 & 284.5 & 218.2 & 231.1 & 216.1 & 472.4 \\
0.6 & 288.3 & 219.7 & 233.2 & 207.1 & 376.1 & 292.4 & 220.1 & 234.0 & 216.3 & 472.2 \\
0.7 & 295.7 & 221.6 & 234.5 & 214.8 & 371.3 & 301.1 & 221.2 & 234.6 & 215.1 & 474.5 \\
0.8 & 302.6 & 236.7 & 246.6 & 257.4 & 342.3 & 308.7 & 232.0 & 243.3 & 257.1 & 398.1 \\
0.9 & 308.8 & 239.7 & 252.7 & 257.1 & 341.6 & 315.3 & 240.1 & 252.8 & 257.1 & 398.3 \\
1.0 & 315.1 & 282.3 & 285.9 & 298.2 & 314.1 & 321.4 & 279.9 & 284.2 & 295.2 & 330.3 \\
\bottomrule
\end{tabular}
%}
\end{table}

\paragraph{Trade-Off Between Cost and Risk Mitigation}
Table~\ref{tab:alpha_lambda} and Figure~\ref{fig:OBJ_lambda} illustrate the impact of increasing $\lambda$ on the objective. As expected, higher $\lambda$ result in a steady increase in the objective across all $\alpha$ levels. This outcome is intuitive since a greater emphasis on extreme scenarios leads to higher overall costs, as the model prioritizes robust solutions that minimize worst-case losses. However, the rate of increase in cost varies with $\alpha$, with higher confidence levels (e.g., $\alpha = 0.90$) showing a steeper rise compared to lower ones (e.g., $\alpha = 0.25$). This suggests that risk-averse settings amplify the effect of risk mitigation efforts, making them increasingly expensive as confidence levels increases.

Conversely, Table~\ref{tab:alpha_lambda} and Figure~\ref{fig:CVAR_lambda} show the diminishing CVaR values as $\lambda$ increases. The reduction in CVaR demonstrates the model’s ability to shift resources toward solutions that mitigate extreme cost realizations. However, beyond a certain point ($\lambda \approx 0.5$), additional increasing $\lambda$ yield only marginal improvements in risk reduction, indicating diminishing returns in CVaR minimization. Notably, this effect is more pronounced for higher $\alpha$ values, reinforcing the idea that extreme risk management strategies become more costly and less effective beyond a specific threshold.

\paragraph{Sensitivity of Objective Value to Confidence Level ($\alpha$)}
Table~\ref{tab:alpha_lambda} and Figure~\ref{fig:OBJ_vs_alpha_fixed_lambda} examine how varying $\alpha$ affects the objective function for fixed $\lambda$ values. Across all $\lambda$ settings, higher $\alpha$ values lead to increased objective values, reflecting the model’s stronger emphasis on extreme costs as confidence levels rise. However, this effect is particularly sensitive to $\lambda$, with higher values ($\lambda = 0.90$) leading to sharp increases in cost as $\alpha$ increases. This behavior confirms that high-risk aversion combined with high confidence levels results in the most conservative and expensive cost allocation strategies, whereas lower $\lambda$ values maintain cost stability across different $\alpha$ levels.

\paragraph{Key Insights}
From the analysis, the following critical insights emerge:
\begin{itemize}
    \setlength{\itemsep}{0pt} % Removes extra space between items
    \setlength{\parskip}{0pt} % Removes extra space between paragraphs within items
    \setlength{\parsep}{0pt}  % Removes extra space between paragraphs outside items
    \item Balancing Cost and Risk: Increasing $\lambda$ improves resilience against extreme cost scenarios but does so at the expense of higher total expenditures. This trade-off is most pronounced for higher $\alpha$ values, where tail-risk mitigation is prioritized.
    \item Diminishing Returns in CVaR Reduction: While higher $\lambda$ values significantly reduce CVaR initially, the rate of reduction slows beyond $\lambda = 0.5$. This suggests that further increases in risk aversion do not necessarily translate into proportionate improvements in extreme risk mitigation.
    \item Growing Cost Disparity Across $\alpha$ Levels: The widening gap in objective values for different $\alpha$ values at higher $\lambda$ levels indicates that high-risk-averse settings amplify sensitivity of the model to changes in confidence levels.
    \item Non-Monotonic CVaR Behavior: At specific configurations, particularly $\lambda = 0.8$ and $\alpha = 0.90$, CVaR exhibits temporary stabilization before further decreasing at $\lambda = 1.0$. This behavior suggests that the model optimally redistributes costs to balance both mean costs and tail risks dynamically.
\end{itemize}

paragraph{Implications for Decision-Makers}
These findings offer valuable insights for strategic decision-making:
\begin{itemize}
    \setlength{\itemsep}{0pt} % Removes extra space between items
    \setlength{\parskip}{0pt} % Removes extra space between paragraphs within items
    \setlength{\parsep}{0pt}  % Removes extra space between paragraphs outside items
    \item For Risk-Averse Stakeholders (e.g., $\alpha = 0.90$): Higher $\lambda$ values effectively minimize CVaR, making them suitable for industries prioritizing financial stability in volatile environments. However, this approach incurs significant additional costs.
    \item For Cost-Sensitive Operators (e.g., $\alpha = 0.25$): Lower $\lambda$ values ensure cost-efficient operations while maintaining moderate exposure to extreme risks, making them ideal for stakeholders with tighter budget constraints.
    \item Optimal Risk Balancing Strategy: The presence of diminishing returns in CVaR reduction suggests that a moderate risk-aversion ($\lambda \approx 0.5$) may offer the best trade-off between cost and worst-case scenario protection.
\end{itemize}

In summary, the results emphasize the importance of balancing cost efficiency with resilience against extreme cost events. The widening gap in objective values across $\alpha$ levels suggests that models incorporating dynamic risk-preference adjustments can better align with diverse operational realities.

\subsubsection{Insights from VSS and EVPI Across Problem Scales}

To assess the benefit of modeling uncertainty, we compute the Value of the Stochastic Solution (VSS) and the Expected Value of Perfect Information (EVPI) across varying problem sizes, defined by combinations of train counts and scenario sets (Table~\ref{tab:vss-cvar}). The table summarizes several key performance metrics: the \textbf{EEV} column reports the expected cost of applying the deterministic plan under uncertainty; \textbf{SS} gives the cost of the optimal stochastic solution; and \textbf{WS} reflects the cost assuming perfect foresight (i.e., full scenario realization). The \textbf{VSS} quantifies the benefit of adopting a stochastic model over a deterministic one, while \textbf{EVPI} measures the hypothetical gain from perfect information. The normalized columns \textbf{VSS\,(\%)} scale these gains relative to the stochastic solution cost. Lastly, \textbf{CVaR\textsubscript{SS}} denotes the Conditional Value-at-Risk, capturing expected cost in the worst-case tail under the stochastic policy.

The results show that VSS values are consistently positive, confirming the advantage of stochastic optimization over deterministic planning. For example, with 4 trains and 6 scenarios, the VSS indicates a 3.8\% cost reduction. This advantage tends to increase with more scenarios, particularly in smaller configurations where deterministic policies exhibit greater vulnerability to uncertainty.

EVPI values, however, are substantially higher than VSS in all settings, e.g., \$52.71k in the (4,8) instance, underscoring the large theoretical value of perfect foresight. The persistent gap between EVPI and VSS reflects the limitations of nonanticipative models, which must commit to first-stage decisions before uncertainty is resolved.

\begin{table}[htbp]
    \centering
    \caption{Impact of Stochastic Optimization and Perfect Information Across Problem Sizes (Costs in \$ $\times$ 10\textsuperscript{3})}
    \label{tab:vss-cvar}
    \begin{tabular}{ccccccccc}
        \hline
        \textbf{Trains} & \textbf{Scen} & \textbf{EEV} & \textbf{SS} & \textbf{WS} &
        \textbf{VSS} & \textbf{EVPI} & \textbf{VSS (\%)} & \textbf{CVaR\textsubscript{SS}} \\
        \hline
        3 & 3 & 44.90 & 43.69 & 14.69 & 1.21 & 29.00 & 2.8 & 31.25 \\
        3 & 6 & 63.16 & 61.50 & 19.09 & 1.66 & 42.41 & 2.7 & 48.10 \\
        3 & 8 & 51.11 & 50.31 & 13.26 & 0.80 & 37.05 & 1.6 & 33.62 \\
        4 & 3 & 50.52 & 49.27 & 19.17 & 1.25 & 30.10 & 2.5 & 32.35 \\
        4 & 6 & 71.86 & 69.20 & 22.82 & 2.66 & 46.38 & 3.8 & 50.30 \\
        4 & 8 & 76.16 & 74.66 & 21.95 & 1.50 & 52.71 & 2.0 & 50.54 \\
        7 & 3 & 89.09 & 88.82 & 33.03 & 0.27 & 55.80 & 0.3 & 56.35 \\
        7 & 6 & 126.82 & 121.90 & 40.84 & 4.91 & 81.06 & 4.0 & 100.60 \\
        7 & 8 & 110.22 & 105.97 & 35.49 & 4.25 & 70.48 & 4.0 & 82.40 \\
        \hline
    \end{tabular}
\end{table}

CVaR values remain proportionally reasonable relative to the SS cost, suggesting that the model effectively controls tail risk without overly inflating expected expenditures. \textbf{Key takeaways:} (i) stochastic optimization yields measurable gains over deterministic plans; (ii) perfect information offers significantly greater cost-saving potential; and (iii) increasing system flexibility (via more trains) tends to reduce the marginal value of stochastic modeling.

Computationally, solving the model across these scales remains tractable for small to moderate instances. The two-stage structure with binary first-stage and integer second-stage decisions, along with CVaR risk terms and inter-scenario coupling, renders the problem NP-hard due to their combinatorial structure and the need to jointly evaluate all scenarios. To assess practical solvability, for the smallest case (3 trains, 3 scenarios), the solve time was under 10 seconds. As the problem scaled to 7 trains and 8 scenarios, solve times increased significantly, with some instances taking several hours to reach optimality. These results demonstrate that the model is tractable for moderate problem sizes, but larger real-world applications may benefit from tailored decomposition or heuristic approaches.

\subsubsection{Analysis of Cost Components}

To better understand the cost structure within our optimization framework, we analyze the aggregated cost components across all scenarios: \textit{Supply Cost}, \textit{Transportation Cost}, \textit{Unmet Demand Penalty}, and \textit{Emissions Penalty}. This breakdown addresses Policy Questions 1, 2, and 6 by illustrating how the model balances economic efficiency, sustainability, and operational constraints under uncertainty.

\paragraph{Cost Drivers and Trade-offs}

\begin{figure}[ht]
    \centering
    \captionsetup{justification=centering} % Center captions and subcaptions
    \includegraphics[width=0.5\textwidth]{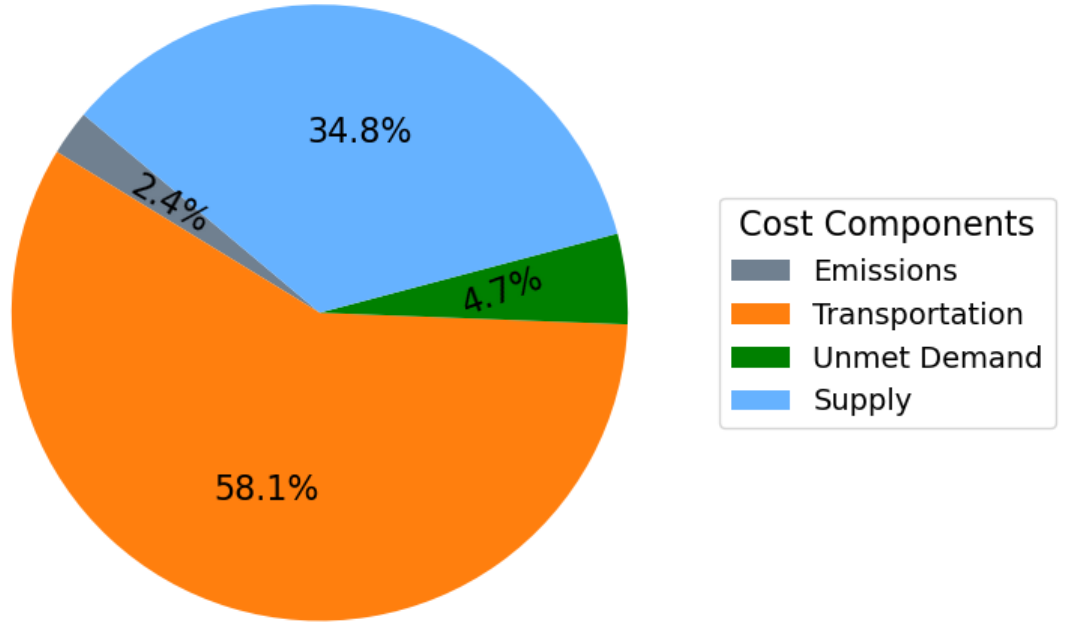}
    \caption{Aggregated cost breakdown across all scenarios.}
    \label{fig:aggregated_costs_pie}
\end{figure}

Figure~\ref{fig:aggregated_costs_pie} shows that \textbf{transportation costs dominate (58.1\%)}, driven by the logistical demands of container movement across the network. Despite the model’s cost-minimization objective, train capacity constraints and network topology elevate these costs. \textbf{Supply costs (34.8\%)} reflect sourcing decisions, constrained by location capacities and their proximity to demand nodes.

\textbf{Unmet demand penalties (4.7\%)} stem from capacity shortfalls, especially under peak demand scenarios. The model prioritizes fulfilling demand where feasible but accepts limited shortfalls to avoid disproportionate cost escalation. \textbf{Emissions penalties (2.4\%)}, while minor, signal the model’s ability to respect sustainability constraints without inflating total costs, though stricter limits could shift this balance.

\paragraph{Insights for Decision-Makers}

The cost breakdown offers several strategic insights:

\begin{itemize}
    \setlength{\itemsep}{0pt} % Removes extra space between items
    \setlength{\parskip}{0pt} % Removes extra space between paragraphs within items
    \setlength{\parsep}{0pt}  % Removes extra space between paragraphs outside items
    \item \textbf{Transportation efficiency is critical}, as it represents the largest share of total costs. Optimizing routes, schedules, and hub operations especially under real-time constraints can yield substantial savings.
    \item \textbf{Sustainability goals can be met at low cost}. The low emissions penalties indicate that the model supports carbon-conscious logistics with minimal financial trade-off, aligning environmental and economic objectives.
    \item \textbf{Capacity constraints drive unmet demand}, but the model shows that smart allocation and flexible sourcing often outperform infrastructure expansion. Adaptive planning and accurate forecasting are key to managing costs under uncertainty.
\end{itemize}

\paragraph{Broader Implications}

This analysis reinforces the value of an integrated planning framework that accounts for operational, economic, and environmental factors. By internalizing emissions and unmet demand into the cost function, the model enables balanced, data-driven decision-making for sustainable intermodal logistics. Future extensions could incorporate dynamic pricing or more granular emissions tracking to further refine trade-off management.

\subsubsection{Impact of Train Capacity Constraints on Demand and Total Cost}

Train capacity plays a critical role in intermodal transport, particularly under uncertain demand and fluctuating spot availability. To explore this, we examine a streamlined setup involving 4 freight trains, 22 intermodal hubs, and 8 demand scenarios, totaling 1013 containers (averaging 32 per train). This configuration captures key uncertainties while allowing us to evaluate how capacity expansion affects cost, demand satisfaction, and system balance.

Table~\ref{tab:capacity_results} summarizes the results. As train capacity increases, unmet demand and total costs decrease sharply up to a point. Expanding capacity from 4 to 8 containers per station yields substantial improvements, but further increases offer only marginal gains. This identifies a capacity threshold beyond which cost savings level off, addressing our sixth policy question on the diminishing returns of capacity investment.

\begin{table}[h]
\scriptsize
    \centering
    \caption{Impact of Train Capacity Constraints on Cost and Demand Satisfaction}
    \label{tab:capacity_results}
    \begin{tabular}{ccccccc}
        \hline
        \textbf{Train Capacity} & \textbf{Total Cost} & \textbf{Unmet D} & \textbf{D Met (\%)} & \textbf{Max Unmet D} & \textbf{Min Unmet D} & \textbf{Std Dev of Unmet D} \\
        \hline
        4  & 249,400  & 453  & 55.28  & 18  & 10  & 3.10  \\
        5  & 214,015  & 293  & 71.08  & 14  & 5   & 2.61 \\
        6  & 182,709 & 178  & 82.43  & 13  & 1   & 3.64 \\
        7  & 156,209  & 70   & 93.09  & 13  & 0   & 2.93 \\
        8  & 137,284  & 13   & 98.72  & 3   & 0   & 0.90  \\
        9  & 127,593  & 9    & 99.11  & 2   & 0   & 0.67 \\
        10 & 121,681  & 9    & 99.11  & 2   & 0   & 0.67 \\
        \hline
    \end{tabular}
\end{table}

\paragraph{System Performance and Trade-offs}  
The performance improvements stem from resolving key bottlenecks:

\begin{itemize}
    \setlength{\itemsep}{0pt} % Removes extra space between items
    \setlength{\parskip}{0pt} % Removes extra space between paragraphs within items
    \setlength{\parsep}{0pt}  % Removes extra space between paragraphs outside items
    \item \textbf{Low-capacity bottlenecks:} When capacity is tight, container allocation is restricted, leading to penalties for unmet demand. Increasing capacity reduces these penalties significantly.
    \item \textbf{Demand volatility:} With 8 demand scenarios, capacity must absorb peaks and troughs. At low levels, trains struggle to meet demand surges; higher capacities offer flexibility and reduce cost variability.
    \item \textbf{Diminishing returns:} Beyond 8 containers, most demand is already met, and additional capacity does little to reduce costs further, as other constraints—such as supply limits and schedule alignment become more influential.
\end{itemize}

\begin{figure}[!htp]
    \centering
    \captionsetup{justification=centering} % Center caption alignment
    \includegraphics[width=0.5\textwidth, height=4.5cm]{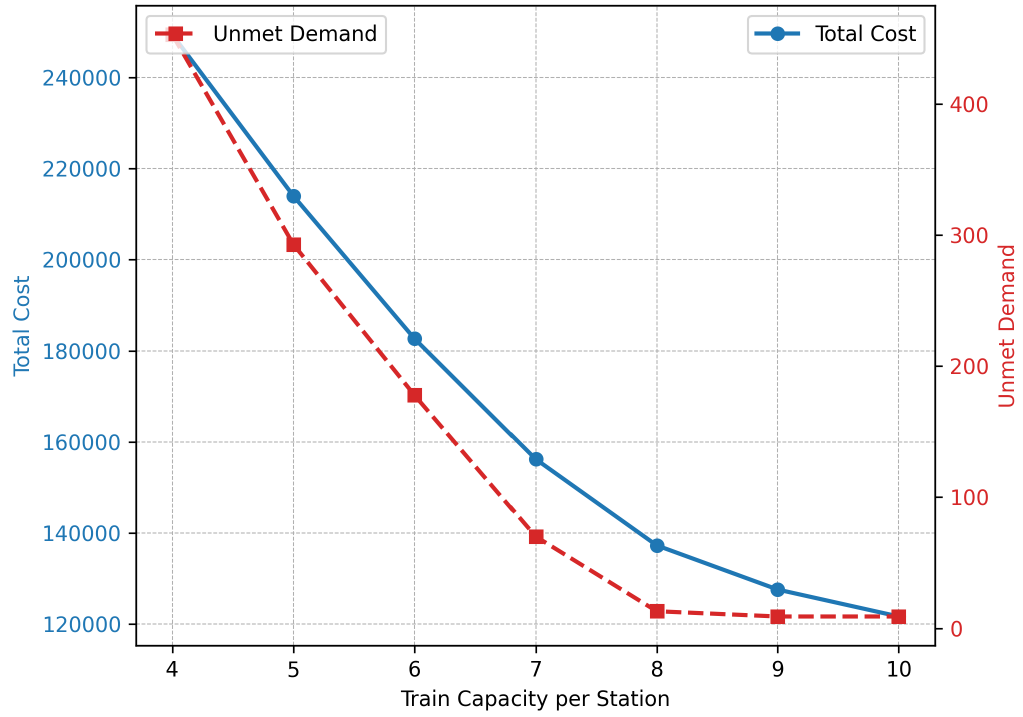}
    \caption{Total Cost and Unmet Demand vs. Train Capacity.}
    \label{fig:cost_vs_unmet}
\end{figure}

\paragraph{Cost vs. Unmet Demand}  
Figure~\ref{fig:cost_vs_unmet} illustrates how cost and unmet demand decline with rising capacity. The steep drop between 4 and 8 containers highlights the benefit of expansion. Beyond that, improvements taper off. This plateau suggests capacity is no longer the main limiter, supply readiness, and hub synchronization are more pronounced.

The distribution of unmet demand also improves. At 4 containers, some trains face up to 18 unmet units, with a high standard deviation. By capacity 10, unmet demand is nearly eliminated and evenly distributed (maximum of 2 units, standard deviation of 0.67), reflecting a more balanced and efficient system.

\paragraph{Key Insights}  
\begin{itemize}
    \setlength{\itemsep}{0pt} % Removes extra space between items
    \setlength{\parskip}{0pt} % Removes extra space between paragraphs within items
    \setlength{\parsep}{0pt}  % Removes extra space between paragraphs outside items
    \item Increasing train capacity significantly improves performance but only up to a point. After 8 containers per station, gains flatten due to other operational limits.
    \item At lower capacities, unmet demand varies widely across trains, indicating congestion and imbalance. Higher capacities smooth out these differences.
    \item The model effectively balances cost minimization with demand satisfaction, optimizing under uncertainty.
    \item Given the unpredictable nature of spot capacity (affected by external demand and network congestion), strategic planning should focus on adaptability and not just expansion.
\end{itemize}

In summary, while capacity expansion enhances flexibility and reliability, it is not a silver bullet. Beyond a certain point, improvements depend more on adaptive scheduling, coordinated resource allocation, and dynamic response to real-time constraints.

\subsubsection{Sensitivity of Cost Metrics to Risk Preferences and Scenario Variability}

\begin{figure}[!htp]
    \centering
    \captionsetup{justification=centering} % Center caption alignment
    \begin{minipage}[t]{0.45\textwidth}
        \centering
        \subfloat[$\alpha = 0.25$\label{fig:OBJ_vs_ETC_alpha25}]{
            \includegraphics[width=\textwidth, height=4cm]{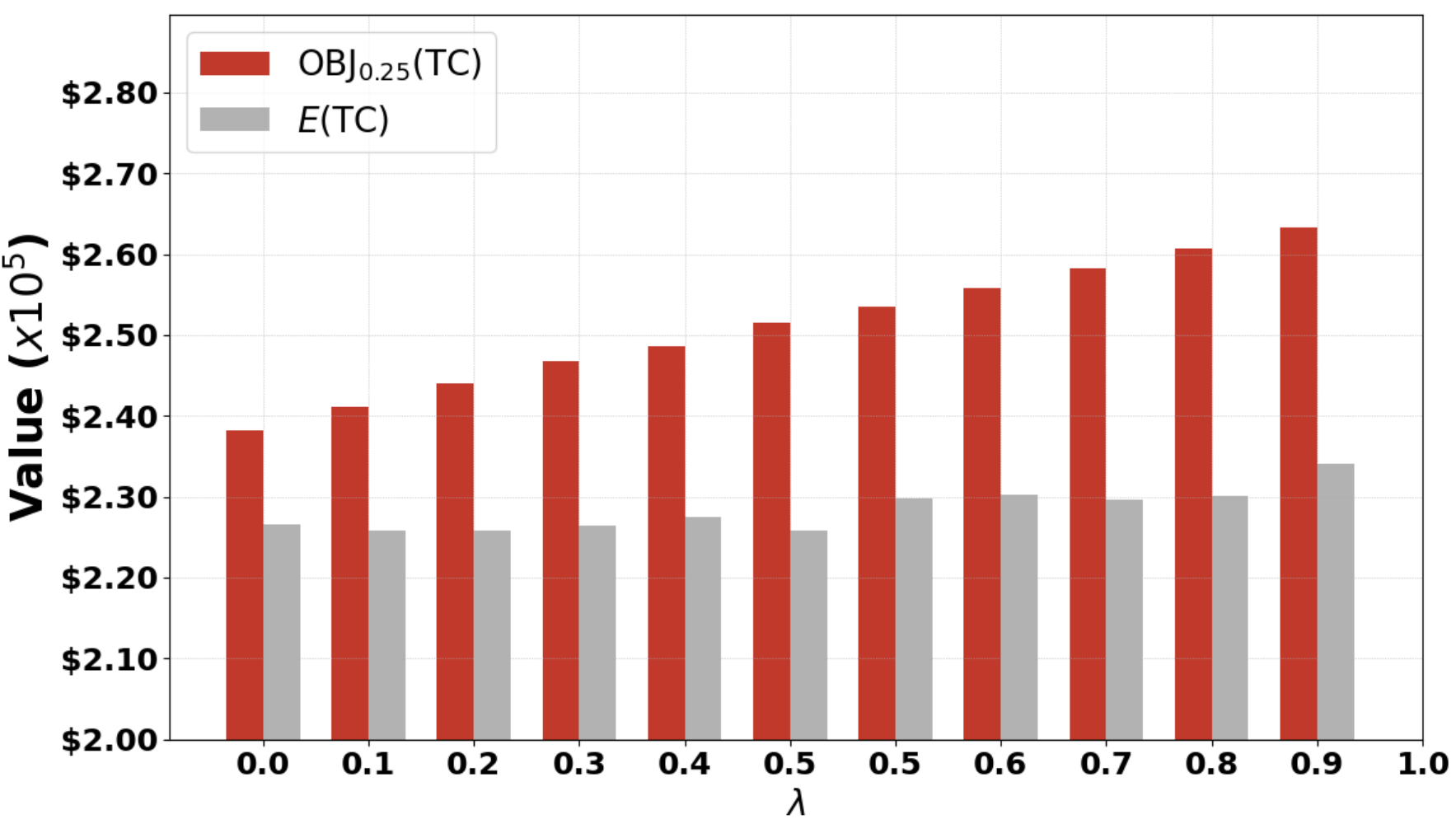}
        } \\
        \vspace{0.3cm} % Add space between figures in the same column
        \subfloat[$\alpha = 0.50$\label{fig:OBJ_vs_ETC_alpha50}]{
            \includegraphics[width=\textwidth, height=4cm]{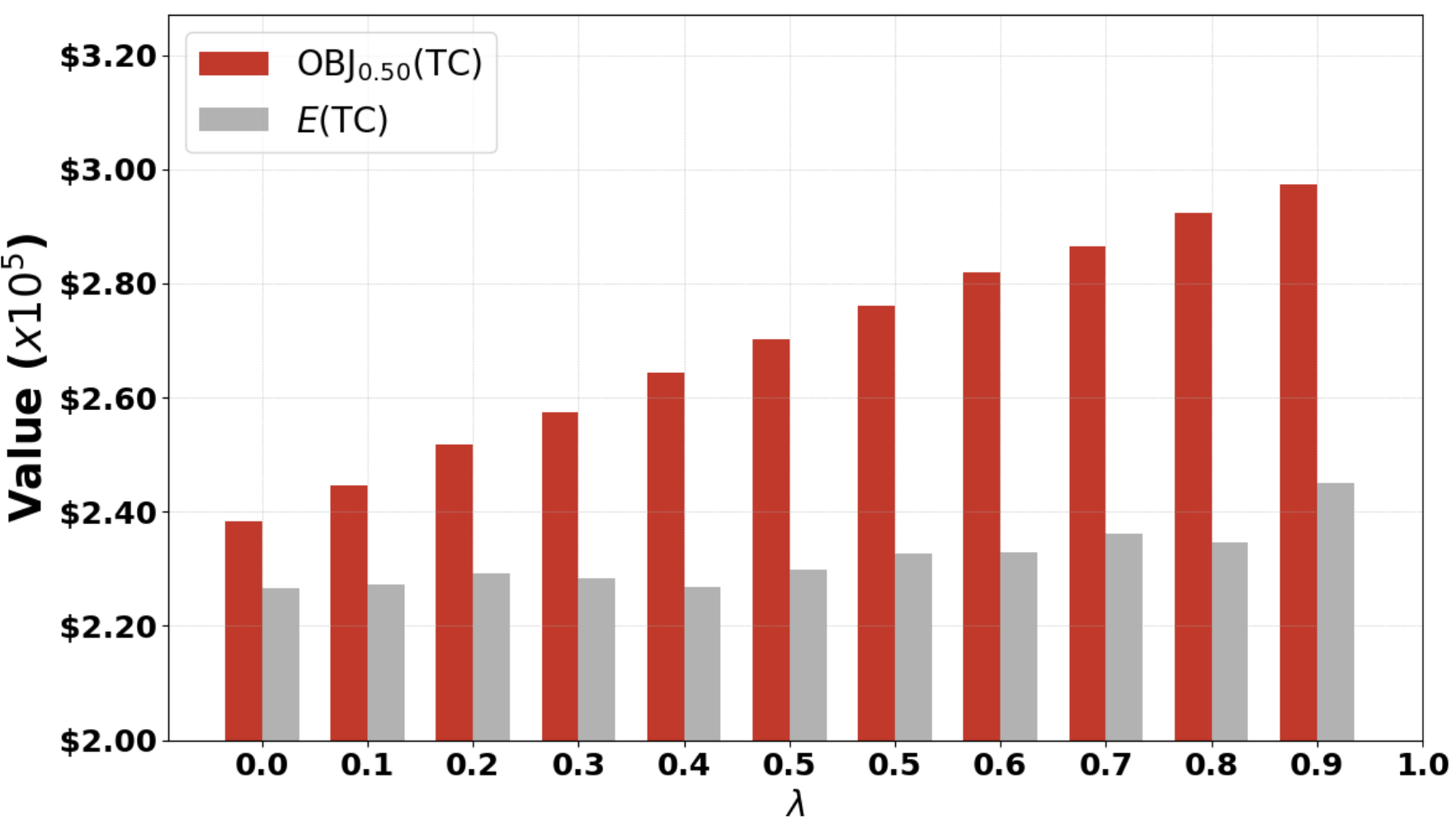}
        }
    \end{minipage}
    \hspace{0.2cm} % Reduce space between columns
    \begin{minipage}[t]{0.45\textwidth}
        \centering
        \subfloat[$\alpha = 0.75$\label{fig:OBJ_vs_ETC_alpha75}]{
            \includegraphics[width=\textwidth, height=4cm]{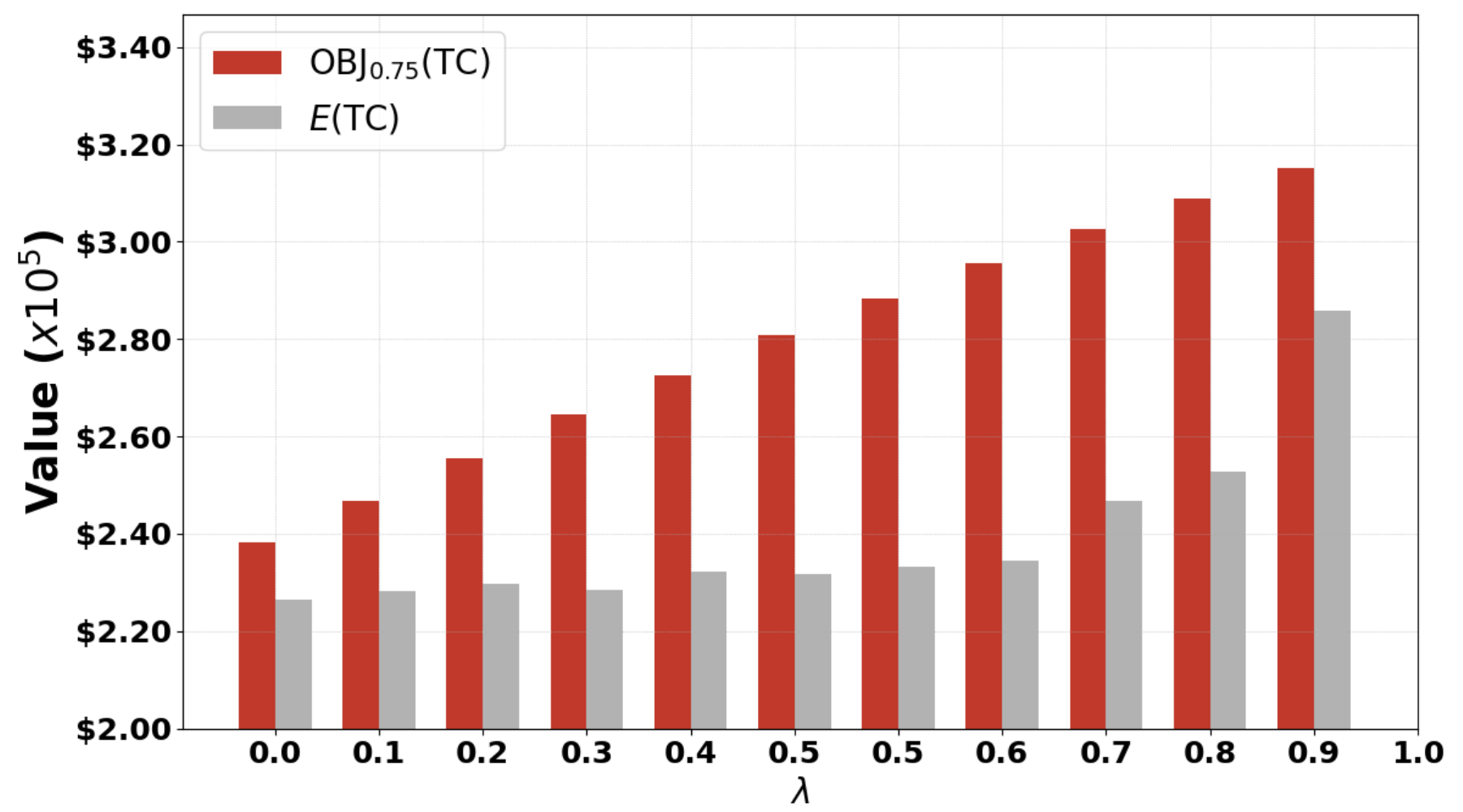}
        } \\
        \vspace{0.3cm} % Add space between figures in the same column
        \subfloat[$\alpha = 0.90$\label{fig:OBJ_vs_ETC_alpha90}]{
            \includegraphics[width=\textwidth, height=4cm]{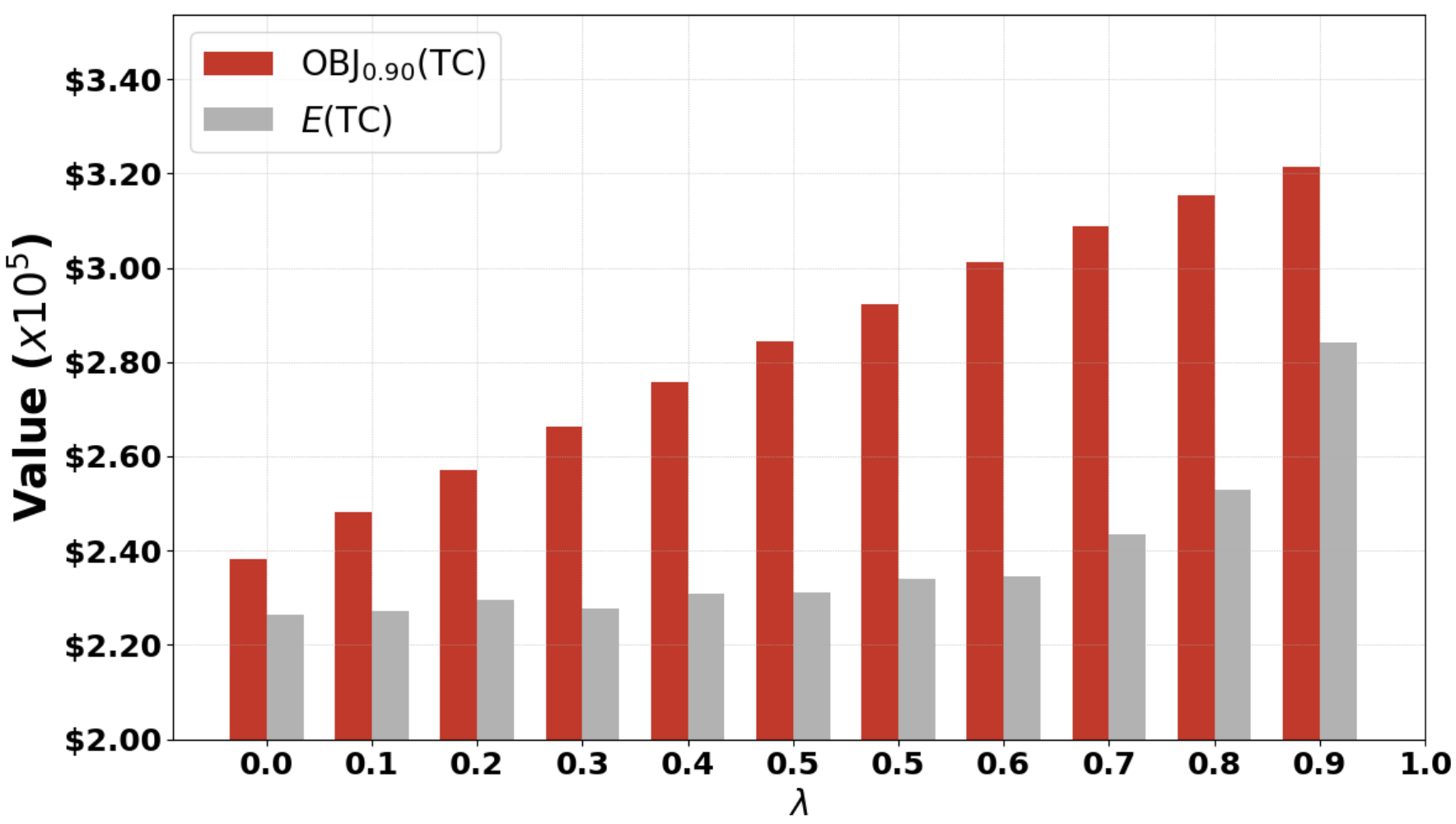}
        }
    \end{minipage}
    \caption{Variation of OBJ Values with E(TC) for Different $\alpha$ Levels.}
    \label{fig:Obj_ETC}
\end{figure}

\begin{figure}[!htp]
    \centering
    \captionsetup{justification=centering} % Center captions and subcaptions
    \subfloat[3D Surface Plot of Total Cost\label{fig:ETC_surface}]{
        \includegraphics[width=0.48\textwidth, height=4.5cm]{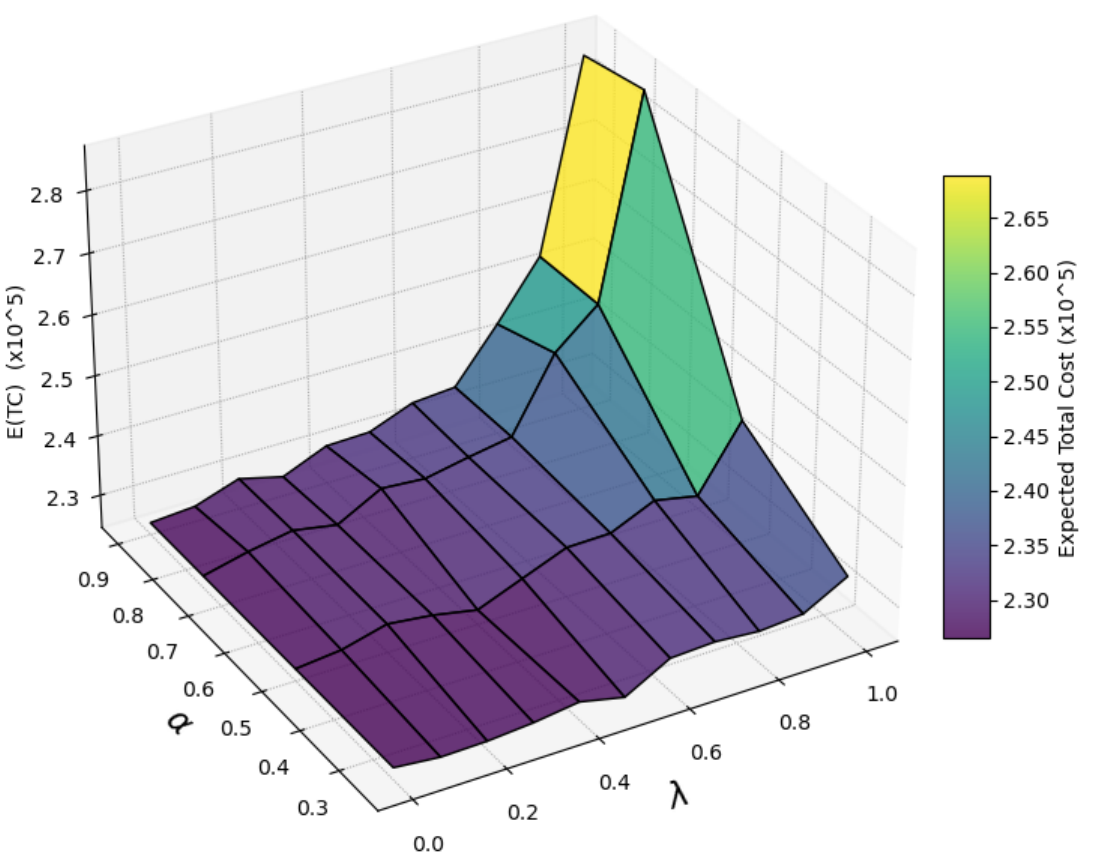}
    }
    \hfill
    \subfloat[Contour Plot of E(TC) Values\label{fig:ETC_contour}]{
        \includegraphics[width=0.48\textwidth, height=4.5cm]{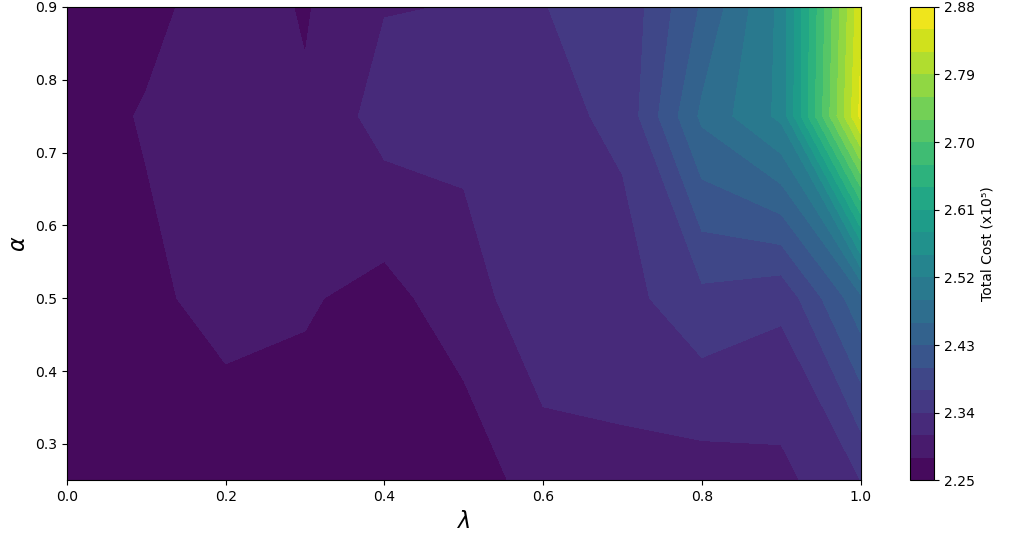}
    }
    \caption{Expected total cost values as a function of $\lambda$ and $\alpha$ illustrating the sensitivity of total cost to changes in $\lambda$ and $\alpha$.}
    \label{fig:ETC_plots}
\end{figure}

A central challenge in stochastic intermodal planning is balancing cost efficiency with robustness across uncertain scenarios. While previous sections explored how $\lambda$ and $\alpha$ shape total cost and CVaR, this subsection focuses on the model’s sensitivity to cost variability and stability. Specifically, it examines the divergence between the objective value (OBJ) and expected total cost $E(TC)$, the distribution of costs across scenarios, and how these factors shift with increasing risk aversion. Figures~\ref{fig:Obj_ETC} and \ref{fig:ETC_plots} visualize these dynamics.

This analysis informs Policy Questions 1, 3, and 7 by illustrating the implications of risk preferences on cost predictability and revealing the diminishing returns of extreme risk aversion.

Figure~\ref{fig:Obj_ETC} shows that as $\lambda$ increases, particularly at higher $\alpha$, the gap between OBJ and $E(TC)$ widens. At low $\lambda$, both metrics remain closely aligned, suggesting that cost minimization dominates the objective function. However, as $\lambda$ grows, CVaR becomes more influential, indicating a shift in the model's focus toward managing worst-case outcomes. This divergence captures the cost-resilience trade-off inherent in conservative planning.

Scenario-level variability, illustrated in Figure~\ref{fig:ETC_plots}, further confirms this trend. At low $\lambda$, the model achieves lower average costs but with significant fluctuation across scenarios. As risk aversion increases, cost volatility declines, producing more stable financial outcomes. Yet, this stabilization levels off around $\lambda \approx 0.8$, beyond which additional risk mitigation yields marginal improvements in variability while significantly raising total cost. This plateau suggests that moderate values of $\lambda$ provide the most efficient balance between stability and cost.

The CVaR sensitivity analysis reinforces these observations. Increasing $\lambda$ reduces CVaR, but the rate of reduction diminishes as $\lambda$ approaches 0.5. The impact of risk aversion is more pronounced at higher $\alpha$ (e.g., $\alpha = 0.90$), where tail-end scenarios carry greater weight. These findings indicate that moderate risk preferences typically $\lambda$ between 0.4 and 0.6 offer strong protection against extreme costs without incurring excessive expense. 

Figures~\ref{fig:ETC_surface} and~\ref{fig:ETC_contour} map $E(TC)$ across the $(\lambda, \alpha)$ space, highlighting how predictable cost structures emerge at intermediate risk settings. At low $\lambda$, the model aggressively allocates costs for efficiency, but at the expense of stability. At high $\lambda$, cost structures flatten but total costs surge, making the trade-offs more pronounced.

From a practical standpoint, results emphasize the need to tailor risk settings to organizational priorities. Low-risk ($\lambda \leq 0.3$) may suit stable environments but carry vulnerability to high-impact scenarios. Balanced settings near $\lambda = 0.5$ offer predictable costs with modest premiums, while highly conservative plans ($\lambda \geq 0.8$) deliver resilience at major cost. This flexibility enables stakeholders to calibrate the model according to their risk tolerance and operational needs.

These findings underscore the importance of aligning model risk settings with strategic objectives. They also highlight potential for exploring how risk preferences interact with emissions penalties, capacity constraints, or real-world disruptions like terminal congestion and supply bottlenecks.

\subsubsection{Analyzing the Influence of Emissions Cap on Model Performance}

\begin{figure}[ht]
    \centering
    \captionsetup{justification=centering} % Center captions and subcaptions
    \begin{minipage}[t]{0.45\textwidth}
        \centering
        \subfloat[$\lambda = 0.25$\label{fig:OBJ_vs_E_.25}]{
            \includegraphics[width=\textwidth, height=4.cm]{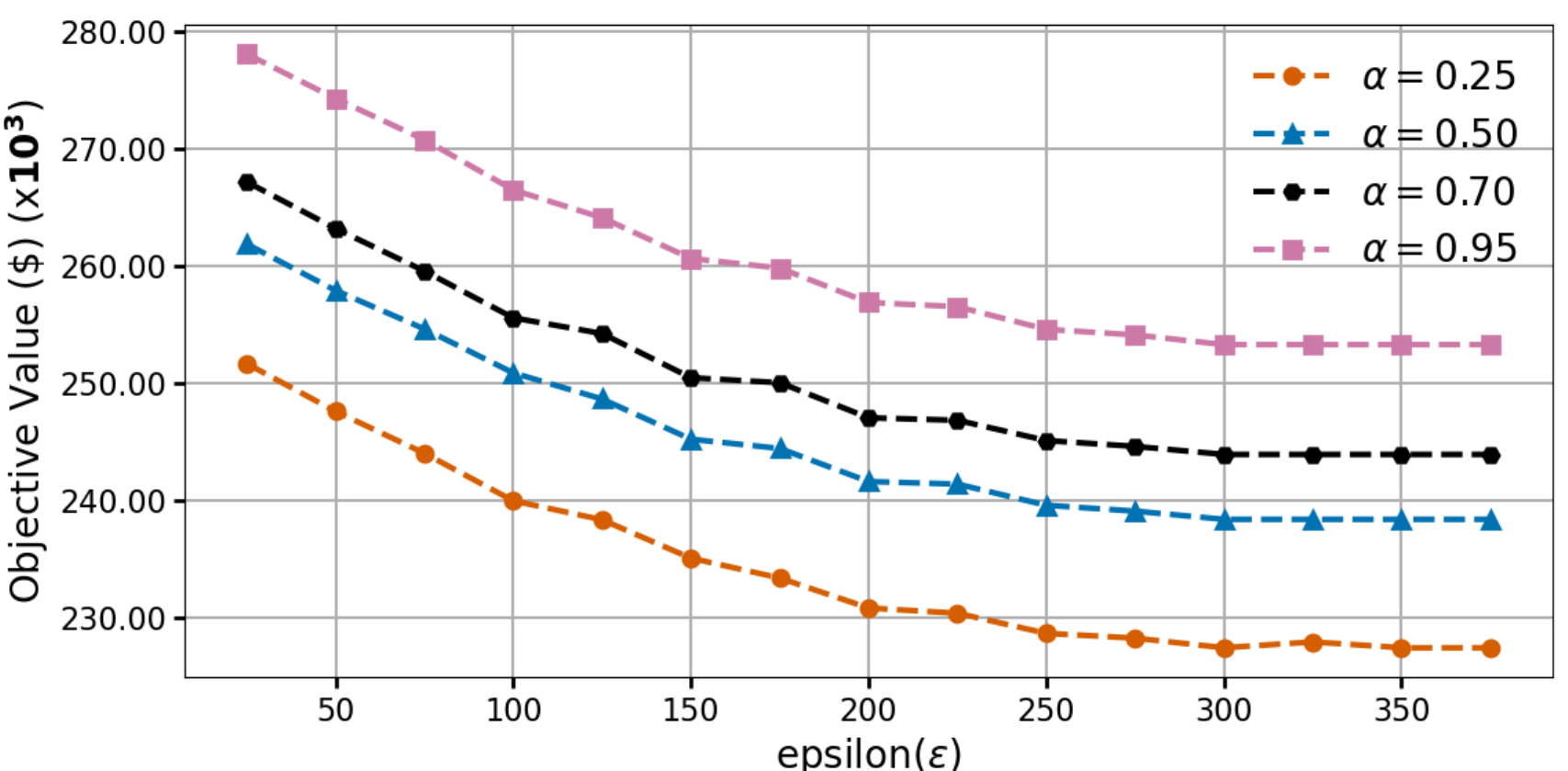}
        } \\
        \vspace{0.3cm} % Add space between figures in the same column
        \subfloat[$\lambda = 0.50$\label{fig:OBJ_vs_E_.5}]{
            \includegraphics[width=\textwidth, height=4.cm]{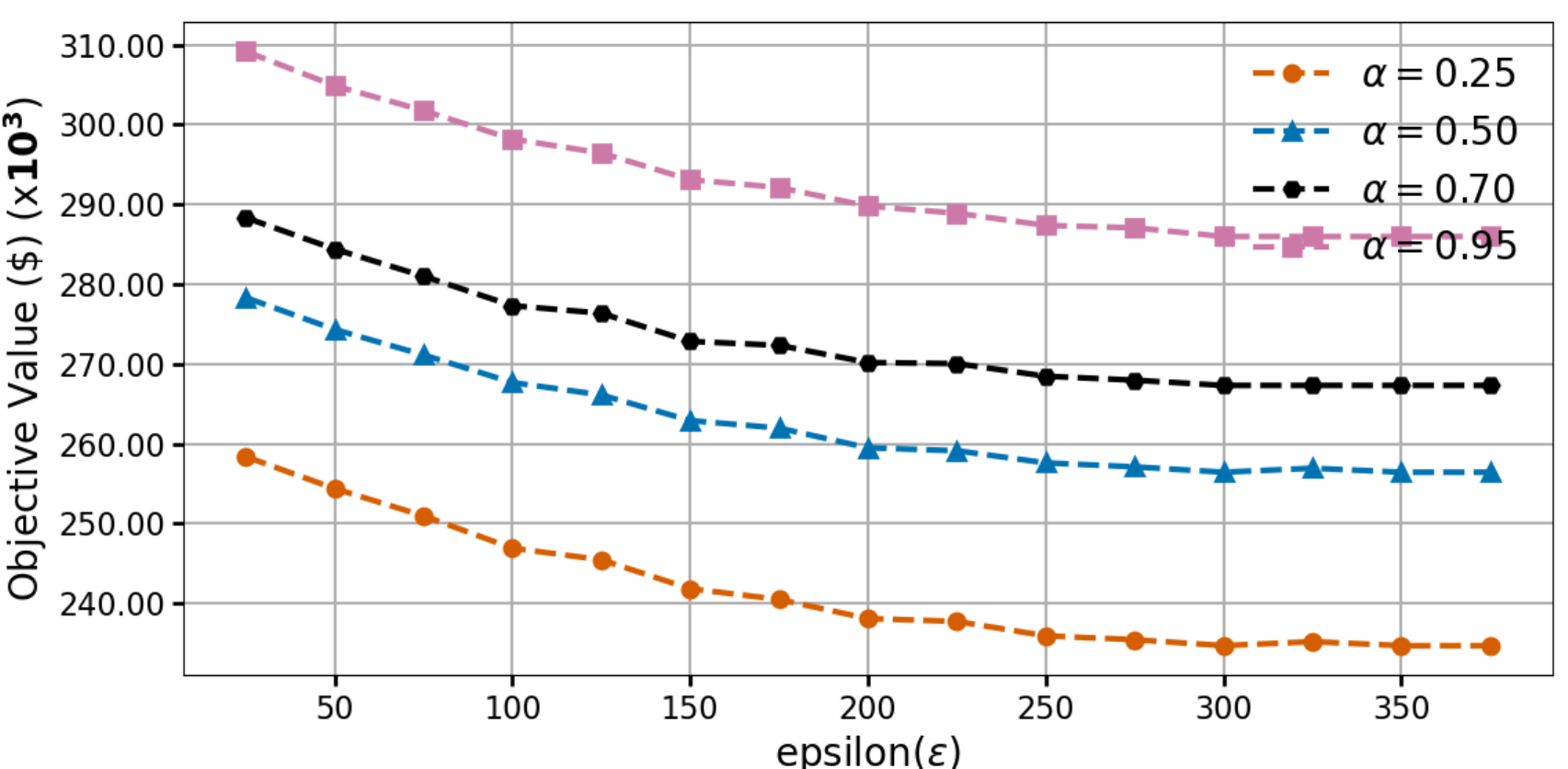}
        }
    \end{minipage}
    \hspace{0.2cm} % Reduce space between columns
    \begin{minipage}[t]{0.45\textwidth}
        \centering
        \subfloat[$\lambda = 0.75$\label{fig:OBJ_vs_E_.75}]{
            \includegraphics[width=\textwidth, height=4.cm]{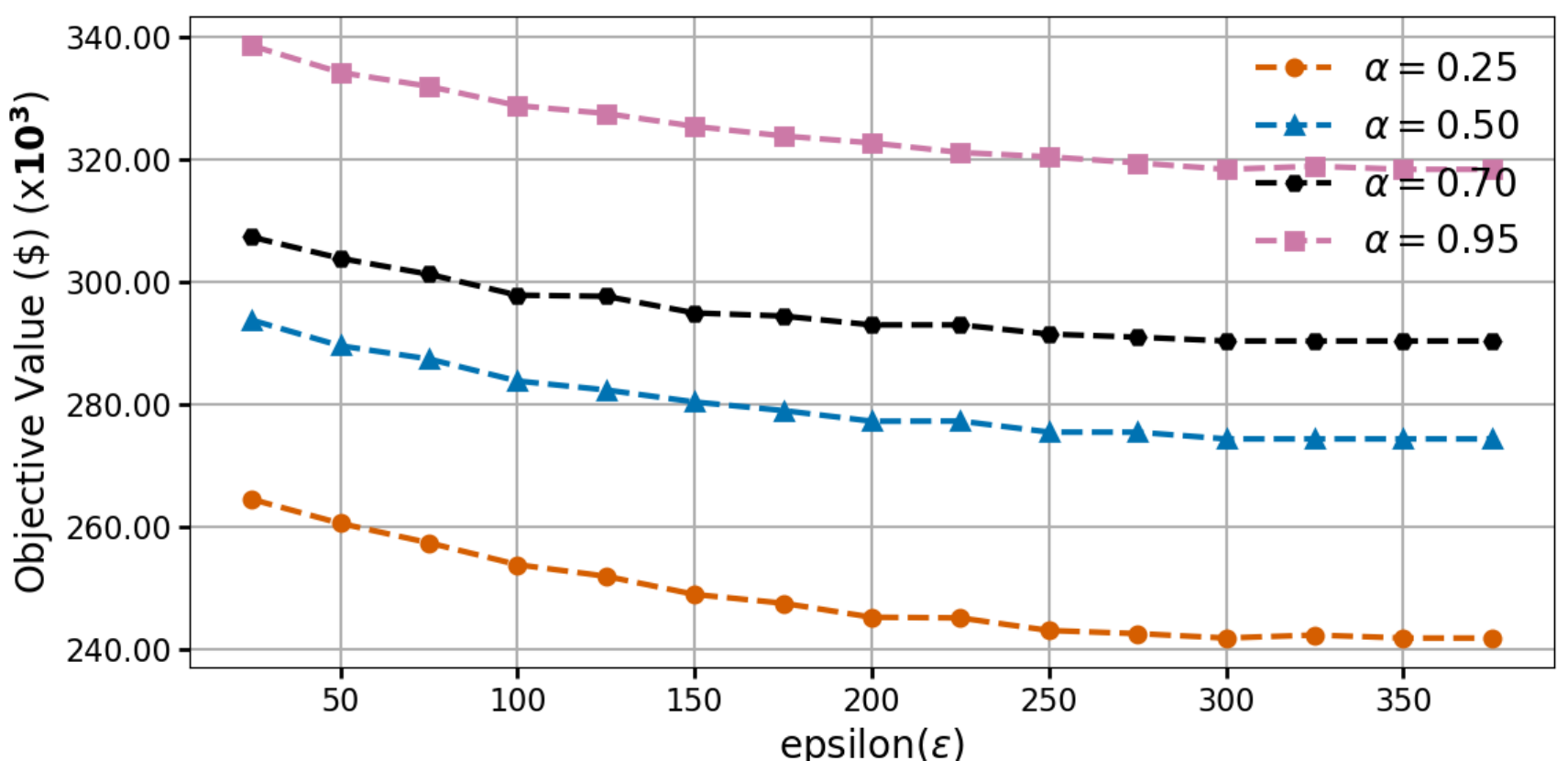}
        } \\
        \vspace{0.3cm} % Add space between figures in the same column
        \subfloat[$\lambda = 0.90$\label{fig:OBJ_vs_E_.9}]{
            \includegraphics[width=\textwidth, height=4.cm]{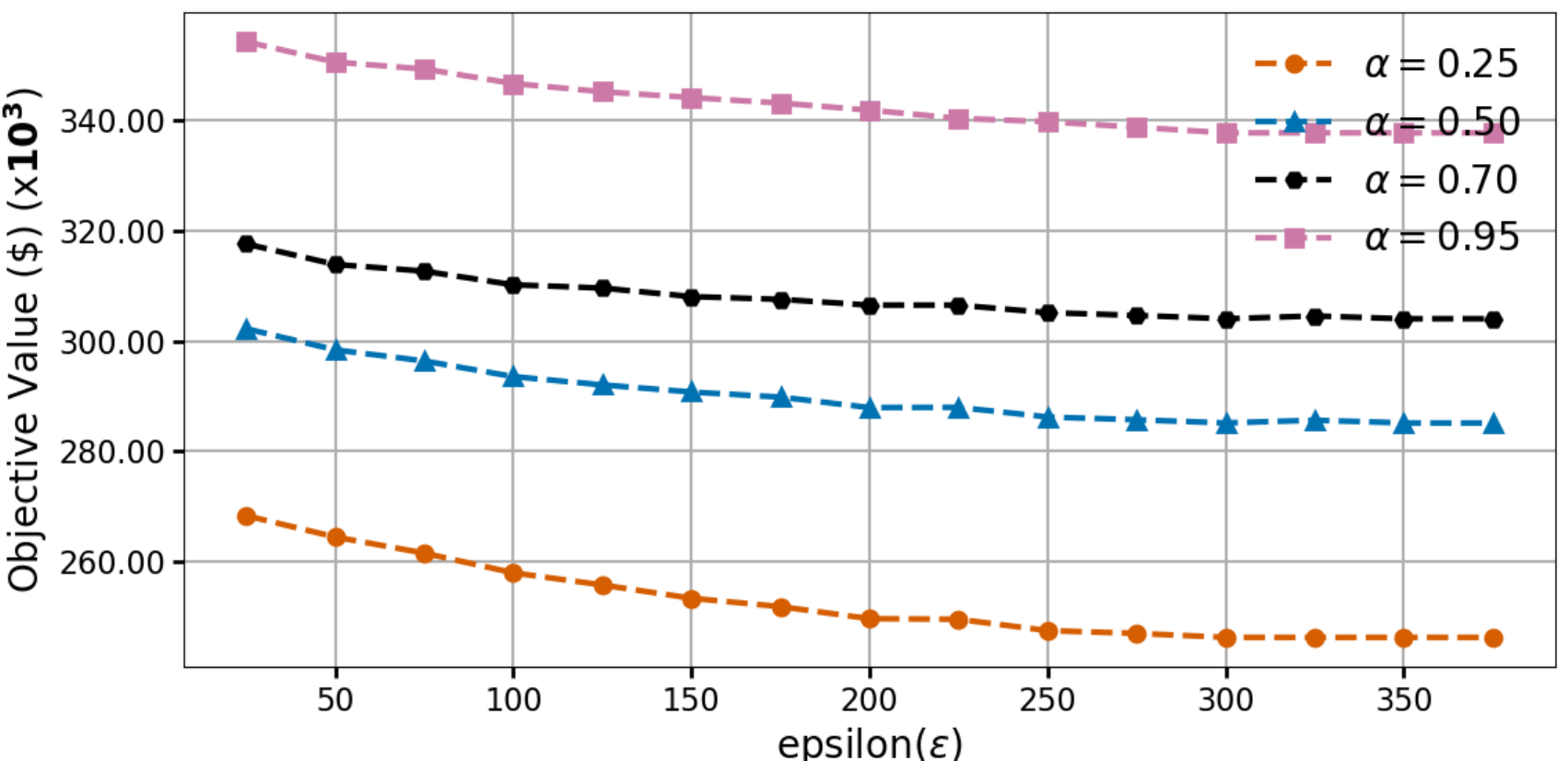}
        }
    \end{minipage}
    \caption{Objective value versus emission cap ($\epsilon$) for different $\alpha$ values at $\lambda$ = 0.25, 0.50, 0.75, and 0.90.}
    \label{fig:emissions_impact}
\end{figure}

Emissions constraints are a critical lever in aligning cost efficiency with sustainability objectives in intermodal transport. This subsection evaluates the impact of emissions caps ($\epsilon$) on total cost under varying risk aversion levels ($\lambda$) and confidence levels ($\alpha$). The results, presented in Table~\ref{tab:objective_values} and visualized in Figures~\ref{fig:OBJ_vs_E_.25} through~\ref{fig:OBJ_vs_E_.9}, reveal how the emissions constraint influences system-wide performance and interacts with the model’s broader cost structure.

As shown in Table~\ref{tab:objective_values}, increasing the $\epsilon$ consistently reduces the objective across combinations of $\lambda$ and $\alpha$. This trend stems from the emissions constraint defined in Equation~\eqref{eq:emissions}, which imposes penalties for exceeding $\epsilon$. Stricter limits (lower $\epsilon$) drive up costs by restricting routing flexibility and prioritizing lower-emission options, while relaxed limits ease those penalties and allow more cost-effective routing. However, the cost benefit of relaxation eventually plateaus, indicating that emissions constraints stop being the primary cost driver once penalties fall below a certain threshold.

\begin{table}[h!]
\centering
\caption{Objective Values for Different Emission Caps with Varying Lambda($\lambda$) and Alpha ($\alpha$) Values (costs in \$ $\times$ 10\textsuperscript{3})}
\label{tab:objective_values}
\resizebox{\textwidth}{!}{
\begin{tabular}{c|cccc|cccc|cccc|cccc}
\toprule
\textbf{} & \multicolumn{4}{c|}{$\lambda = 0.25$} & \multicolumn{4}{c|}{$\lambda = 0.50$} & \multicolumn{4}{c|}{$\lambda = 0.75$} & \multicolumn{4}{c}{$\lambda = 0.90$} \\
\midrule
$\boldsymbol{\epsilon}$ & $\alpha = 0.25$ & $\alpha = 0.50$ & $\alpha = 0.70$ & $\alpha = 0.95$ & $\alpha = 0.25$ & $\alpha = 0.50$ & $\alpha = 0.70$ & $\alpha = 0.95$ & $\alpha = 0.25$ & $\alpha = 0.50$ & $\alpha = 0.70$ & $\alpha = 0.95$ & $\alpha = 0.25$ & $\alpha = 0.50$ & $\alpha = 0.70$ & $\alpha = 0.95$ \\
\midrule
25 & 251.7 & 261.9 & 267.2 & 278.1 & 258.4 & 278.3 & 288.3 & 309.2 & 264.5 & 293.8 & 307.3 & 338.6 & 268.3 & 302.3 & 317.6 & 354.4 \\
50 & 247.6 & 257.9 & 263.2 & 274.3 & 254.4 & 274.3 & 284.4 & 304.9 & 260.6 & 289.6 & 303.9 & 334.2 & 264.4 & 298.4 & 313.9 & 350.7 \\
75 & 244.0 & 254.6 & 259.6 & 270.8 & 251.0 & 271.1 & 281.0 & 301.8 & 257.4 & 287.4 & 301.3 & 331.9 & 261.5 & 296.4 & 312.7 & 349.4 \\
100 & 240.0 & 250.9 & 255.6 & 266.5 & 246.9 & 267.7 & 277.3 & 298.2 & 253.9 & 283.8 & 297.8 & 328.8 & 257.9 & 293.6 & 310.2 & 346.8 \\
125 & 238.4 & 248.7 & 254.3 & 264.1 & 245.4 & 266.1 & 276.4 & 296.4 & 252.0 & 282.4 & 297.7 & 327.5 & 255.7 & 292.0 & 309.6 & 345.3 \\
150 & 235.1 & 245.2 & 250.5 & 260.7 & 241.8 & 262.9 & 272.8 & 293.1 & 249.0 & 280.4 & 294.9 & 325.4 & 253.3 & 290.8 & 308.1 & 344.2 \\
175 & 233.4 & 244.5 & 250.1 & 259.8 & 240.5 & 262.0 & 272.3 & 292.1 & 247.6 & 279.0 & 294.4 & 323.8 & 251.8 & 289.8 & 307.6 & 343.2 \\
200 & 230.9 & 241.6 & 247.1 & 256.9 & 238.1 & 259.5 & 270.2 & 289.8 & 245.3 & 277.3 & 293.0 & 322.7 & 249.6 & 287.9 & 306.5 & 341.9 \\
225 & 230.4 & 241.4 & 246.9 & 256.5 & 237.7 & 259.1 & 270.0 & 288.9 & 245.2 & 277.3 & 293.0 & 321.1 & 249.5 & 288.0 & 306.5 & 340.4 \\
250 & 228.7 & 239.6 & 245.1 & 254.6 & 235.9 & 257.6 & 268.5 & 287.4 & 243.1 & 275.5 & 291.5 & 320.4 & 247.4 & 286.2 & 305.2 & 339.8 \\
275 & 228.3 & 239.1 & 244.6 & 254.1 & 235.4 & 257.1 & 268.0 & 287.1 & 242.6 & 275.5 & 291.0 & 319.4 & 246.9 & 285.7 & 304.7 & 338.8 \\
300 & 227.5 & 238.4 & 243.9 & 253.3 & 234.7 & 256.4 & 267.3 & 286.0 & 241.9 & 274.4 & 290.4 & 318.4 & 246.2 & 285.1 & 304.1 & 337.8 \\
325 & 228.0 & 238.4 & 243.9 & 253.3 & 235.2 & 256.9 & 267.3 & 286.0 & 242.4 & 274.4 & 290.4 & 318.9 & 246.2 & 285.6 & 304.6 & 337.8 \\
350 & 227.5 & 238.4 & 243.9 & 253.3 & 234.7 & 256.4 & 267.3 & 286.0 & 241.9 & 274.4 & 290.4 & 318.4 & 246.2 & 285.1 & 304.1 & 337.8 \\
375 & 227.5 & 238.4 & 243.9 & 253.3 & 234.7 & 256.4 & 267.3 & 286.0 & 241.9 & 274.4 & 290.4 & 318.4 & 246.2 & 285.1 & 304.1 & 337.8 \\
\bottomrule
\end{tabular}
}
\end{table}

The interaction between emissions constraints and risk preferences is also notable. At higher $\lambda$ values, the model places greater emphasis on mitigating extreme costs, making it sensitive to changes in $\epsilon$. This effect intensifies with higher $\alpha$, which shifts focus toward tail-end scenarios. As illustrated in Figures~\ref{fig:OBJ_vs_E_.25}–\ref{fig:OBJ_vs_E_.9}, the cost gap between different emissions caps widens with increasing risk aversion. This demonstrates a clear trade-off: tighter emissions limits can enhance system resilience but at the cost of elevated expenditures, particularly in more conservative planning regimes.

Interestingly, the objective values begin to plateau at higher emissions thresholds, and in some cases, exhibit minor non-monotonic fluctuations. These irregularities shows solver behavior near optimality, where multiple feasible solutions deliver similar outcomes but slightly redistribute cost components such as supply or transport. Once emissions penalties become negligible, the model’s focus naturally shifts toward optimizing other drivers.

These findings confirm that substantial emissions reductions can be achieved with minimal cost trade-offs, particularly when emissions constraints are calibrated effectively. At moderate $\lambda$ and $\alpha$ values, the model jointly balances risk, cost, and sustainability objectives, supporting operational decisions that align with environmental policy. Notably, the results indicate a threshold effect: beyond an emissions cap of approximately 200–250 tons, further relaxation yields little marginal cost benefit. This suggests that carbon caps or pricing schemes targeting this range can serve as effective policy levers, encouraging mode shifts without imposing excessive cost burdens. These insights offer practical value for regulators or enterprises evaluating carbon tax thresholds or emissions allowances, and they support future extensions involving explicit carbon pricing, credit trading mechanisms, or adaptive policy rules.

\subsubsection{Incorporating Intermodal Transfer Time and Cost}

To enhance the applicability and realism of the proposed model, we extend the formulation to include \textit{intermodal transfer time} and \textit{transfer cost}, capturing the handling delays and financial expenditures incurred when transitioning goods from road to rail at intermodal hubs. While logistics companies cannot directly control these factors because they are managed by terminal operators and rail service providers their inclusion provides a clearer understanding of their impact on cost structure and scheduling decisions. Companies must ensure containers arrive at hubs with sufficient lead time to accommodate transfer activities and scheduled train departures.

Since the model is already designed to ensure timely container arrivals at intermodal hubs, additional experiments focusing solely on transfer components will not be carried out in this study. The key distinction in this extension lies in the addition of transfer costs per container and the time required for modal shifts. Although stakeholders cannot directly optimize intermodal transfers, incorporating these parameters strengthens the model’s applicability. Future research could build on this extension by modeling the entire transportation process from the supply location to final consumer, offering deeper insights into system-wide coordination across multiple hubs and transportation modes.

\paragraph{Motivation and New Parameters}
Intermodal transfer activities such as unloading, sorting, and reloading introduce delays and costs that influence operational efficiency and scheduling. While not directly controllable by stakeholders in this study, these parameters provide critical insights into cost interplay and constraints within intermodal logistics. Their inclusion enhances the model’s predictive accuracy and practical applicability.

The following parameters are introduced to represent the effects of intermodal transfers:  
\begin{itemize}
    \item $\tau_{ij}^{\text{tr}}$: Transfer time required to move goods from road to rail at intermodal station $j$, originating from location $i$. This parameter reflects delays from handling, processing, and equipment usage at the hub.  
    \item $C_{ij}^{\text{tr}}$: Transfer cost incurred for moving goods from road to rail at station $j$, originating from location $i$. This cost captures labor, equipment fees, and infrastructure charges associated with the intermodal transfer.
\end{itemize}

\paragraph{Enhanced Objective Function.}
The inclusion of intermodal transfer cost $C_{\text{tr}}$ enhances the original objective function by explicitly capturing these additional costs. The revised objective function is as follows:  
\vspace{-1.5em}

\begin{align}
\text{Min} \quad Z =  & \sum_{i \in I} O_i y_i + (1 - \lambda) \sum_{\omega \in \Omega} p_\omega \left( \sum_{i \in I} \sum_{j \in J} \sum_{n \in N} \sum_{t \in T} (C_{ij} + C^{\text{tr}}_{ij}) x_{ijn\omega t} + \pi_\omega \sum_{n \in N} U_{n\omega} \right) \nonumber \\
& + \lambda \text{CVaR} + \sum_{t \in T} \rho \eta_\omega \label{eq:objective_updated}
\end{align}

where $C_{ij}^{\text{tr}}$ accounts for the cost incurred during modal transitions at intermodal hubs. This modification ensures that transfer costs are evaluated alongside other operational components, allowing the model to account for all cost factors within intermodal logistics systems.

\paragraph{Updated Travel Time Constraints.}
The introduction of transfer time $\tau_{ij}^{\text{tr}}$ requires adjustments to the travel time constraints. Goods transported from origin $i$ to station $j$ must now respect the total travel time, which includes both transit time $\tau(i, j)$ and transfer time $\tau_{ij}^{\text{tr}}$. The updated constraint is as follows:  
\vspace{-1.5em}

\begin{align}
    \sum_{i \in I} x_{ijn \omega t} \leq y_i, \quad \forall j \in J, \, \forall n \in N, \, \forall \omega \in \Omega, \, \forall t \in T, \text{ if } t \geq \tau_{ij} + \tau_{ij}^{\text{tr}}.
\label{eq:transfer_time_constraint}
\end{align}

Additionally, goods cannot be transferred to the station earlier than the cumulative travel time allows, as enforced by the constraint below:

\vspace{-1.5em}
\begin{align}
    \sum_{i \in I} x_{ijn \omega t} = 0, \quad \forall j \in J, \, \forall n \in N, \, \forall \omega \in \Omega, \, \forall t \in T, \text{ if } t < \tau_{ij} + \tau_{ij}^{\text{tr}}.
\end{align}
These constraints ensure that travel time, including intermodal transfers, aligns with trains schedules.

\paragraph{Stakeholder Implications and Future Extensions.}
For logistics and manufacturing companies, intermodal transfer costs and delays are external factors determined by rail operators and hub administrators. While stakeholders cannot directly optimize these components, ensuring timely container arrivals at intermodal hubs allows the model to incorporate these constraints into the planning process effectively. Including intermodal transfer parameters lays the groundwork for future research on system-wide logistics optimization. A more comprehensive top-down approach could explicitly model transfer activities, capturing their interdependencies and impact across the broader supply chain.

In summary, while intermodal transfer costs and time remain fixed external inputs in the current model, their integration enhances the model’s realism and applicability. This extension not only strengthens its alignment with real-world constraints but also sets the stage for further investigations into intermodal coordination and optimization.

\section{Conclusion and Future Work} \label{Conclusion}

This study introduces a novel two-stage stochastic optimization model for road-rail intermodal freight transportation, by integrating demand and capacity uncertainty, risk management, and sustainability objectives in a unified framework. Our model uniquely combines aspects of the intermodal network such as shipment allocation and carbon-efficient routing, offering a comprehensive approach to optimizing intermodal logistics under uncertainty. One of the key contributions of this study is the incorporation of real-time adaptive responses to stochastic elements such as train capacities and demand fluctuations while dynamically adjusting container allocation based on realized demand and capacity constraints. This adaptive approach ensures that decisions are aligned with operational realities, providing a level of flexibility and responsiveness that has not been previously achieved in intermodal freight optimization.

Another distinctive feature of our model is its dual emphasis on cost efficiency and environmental sustainability. Our model embeds carbon emission constraints directly into the optimization framework. By penalizing excess emissions and incorporating these penalties into the objective function, we ensure that carbon footprint reduction is a central consideration in the decision-making process. This approach aligns with regulatory frameworks, offering a practical solution for stakeholders seeking to comply with environmental standards. Furthermore, our model introduces a risk-averse optimization framework that leverages CVaR to manage extreme cost scenarios. While CVaR has been widely used in financial and supply chain contexts, its application to intermodal freight transportation under demand and capacity uncertainty is novel. Our findings demonstrate that the integration of CVaR not only enhances system resilience but also provides decision-makers with the flexibility to balance cost efficiency with risk aversion.

Building on these contributions, future work could extend the model to incorporate dynamic routing and variable travel time. In this study, we assumed fixed road transport times for tractability and consistency with short-haul, repetitive intermodal routes often assumed to be fixed \citep{basallo2021planning}. The intermodal network structure and routing paths were also predetermined. While these assumptions are suitable for strategic planning, they limit responsiveness to real-time disruptions such as congestion or weather events. Acknowledging these limitations, future work could explore real-time optimization or simulation-based methods to enhance the model’s adaptability to dynamic conditions. Additionally, since we used a commercial solver, exploring heuristic or decomposition-based approaches may improve scalability for large real-world instances. Second, future work could investigate integration of predictive analytics and machine learning to improve demand forecasting and capacity planning. By leveraging real-time data and predictive models, the optimization framework could be further enhanced to anticipate and respond to fluctuations in demand and capacity more effectively. Finally, the model could be expanded to include multi-objective optimization techniques that simultaneously consider additional sustainability metrics, such as energy consumption. This would provide a more holistic approach to sustainable freight transport, addressing a wider range of environmental concerns.

In conclusion, this study advances the field of intermodal freight transportation by addressing critical gaps in the literature and offering a robust, flexible, and sustainable optimization framework. The unique findings and contributions of this work provide a foundation for further exploration and innovation in the pursuit of efficient, resilient, and environmentally sustainable freight transportation systems.

\section*{Declaration of interests}
The authors have no competing interests to declare that are relevant to the content of this article.

\section*{Acknowledgment}
\label{sec:acknowledgment}
This work was partially supported by the U.S. Department of Energy's Advanced Research Projects Agency-Energy (ARPA-E) under the project (\#DE-AR0001780) titled \textit{``A Cognitive Freight Transportation Digital Twin for Resiliency and Emission Control Through Optimizing Intermodal Logistics"}.

\end{document}